\documentclass{amsart}
\usepackage{amsfonts,amssymb,amscd,amsmath,enumerate,verbatim}
\usepackage[all]{xy}
\SelectTips{eu}{}
\xyoption{curve}
\CompileMatrices
\usepackage[colorlinks]{hyperref}

\newcommand{\xla}{\xleftarrow}
\newcommand{\lra}{\longrightarrow}

\newcommand{\xra}{\xrightarrow}

\newcommand{\card}{{\operatorname{card}}}
\newcommand{\even}{{\operatorname{even}}}
\newcommand{\odd}{{\operatorname{odd}}}
\newcommand{\op}{{\operatorname{op}}}

\newcommand{\la}{\langle}
\newcommand{\ra}{\rangle}
\newcommand{\ov}{\overline}

\newcommand{\wh}{\widehat}
\newcommand{\wt}{\widetilde}

\newcommand{\les}{{\scriptscriptstyle\leqslant}}
\newcommand{\ges}{{\scriptscriptstyle\geqslant}}

\newcommand{\vge}{{\scriptstyle{\mathbf\backslash}\!\bigvee}}
\newcommand{\vgu}{{\bigcup\!\scriptstyle{\mathbf\vert}}}

\newcommand{\dcat}[1]{{\mathsf{D}(#1)}}
\newcommand{\lotimes}[1]{\otimes^{\bf L}_{#1}}

\newcommand{\Rhom}[3]{{\bf R}\!\Hom{#1}{#2}{#3}}

\newcommand{\shift}{{\sf\Sigma}}
\newcommand\col{\colon}
\newcommand{\cone}[1]{{\mathsf{cone}}(#1)}
\newcommand\dd{\partial}
\newcommand{\hh}[1]{\operatorname{H}(#1)}
\newcommand{\HH}[2]{\operatorname{H}_{#1}(#2)}

\newcommand\ZZ{\operatorname{Z}}

\newcommand{\Tor}[4]{\operatorname{Tor}_{#1}^{#2}(#3,#4)}
\newcommand{\Ext}[4]{\operatorname{Ext}^{#1}_{#2}(#3,#4)}
\newcommand{\Hom}[3]{\operatorname{Hom}_{#1}(#2,#3)}

\newcommand\Ker{\operatorname{Ker}}

\newcommand\Image{\operatorname{Im}}

\newcommand{\Ann}{\operatorname{Ann}}

\newcommand\codim{\operatorname{codim}}

\newcommand\edim{\operatorname{edim}}
\newcommand\depth{\operatorname{depth}}
\newcommand\height{\operatorname{height}}
\newcommand{\indx}[1]{\operatorname{index}#1}

\newcommand{\reg}[2]{\operatorname{reg}_{#1}#2}
\newcommand\supheight{\operatorname{super\,height}}
\newcommand\length{\operatorname{\ell}}
\newcommand{\lol}[2]{{\operatorname{\ell\hspace{-1pt}\ell}}_{#1}#2}
\newcommand{\ord}{\operatorname{ord}}

\newcommand{\rank}{\operatorname{rank}}
\newcommand{\frank}{\operatorname{f-rank}}
\newcommand{\cnr}{\operatorname{cf-rank}}

\newcommand\idmap{\operatorname{id}}

\newcommand\pd{\operatorname{proj\,dim}}

\newcommand\gldim{\operatorname{gl\,dim}}

\newcommand\gr{\operatorname{gr}}

\newcommand\fm{{\mathfrak m}}
\newcommand\fn{{\mathfrak n}}

\newcommand\fp{{\mathfrak p}}
\newcommand\fq{{\mathfrak q}}

\newcommand\bsa{{\boldsymbol a}}

\newcommand\bsg{{\boldsymbol g}}

\newcommand\bsx{{\boldsymbol x}}
\newcommand\bsy{{\boldsymbol y}}

\newcommand\BF{{\mathbb F}}
\newcommand\BN{{\mathbb N}}

\newcommand\BQ{{\mathbb Q}}

\newcommand\BZ{{\mathbb Z}}

\newcommand\ca{{\mathsf A}}
\newcommand\cb{{\mathsf B}}
\newcommand\cc{{\mathsf C}}
\newcommand\ce{{\mathsf E}}
\newcommand\cf{{\mathsf F}}

\newcommand\cs{{\mathsf S}}
\newcommand\ct{{\mathsf T}}
\newcommand\cu{{\mathsf U}}

\newcommand\SF{{\mathsf f}}
\newcommand\SG{{\mathsf g}}
\newcommand\SH{{\mathsf h}}

\newcommand\SII{{\mathsf i}}
\newcommand\SJ{{\mathsf j}}
\newcommand\SK{{\mathsf k}}

\newcommand\ST{{\mathsf t}}
\newcommand{\eps}{{\varepsilon}}
\newcommand{\vf}{{\varphi}}
\newcommand{\vk}{{\varkappa}}
\newcommand{\var}{{\hskip1pt\vert\hskip1pt}}

\theoremstyle{plain}
\newtheorem{theorem}{Theorem}[section]

\newtheorem{proposition}[theorem]{Proposition}
\newtheorem{lemma}[theorem]{Lemma}

\newtheorem{corollary}[theorem]{Corollary}

\newtheorem{subtheorem}{Theorem}[theorem]

\newtheorem{sublemma}[subtheorem]{Lemma}

\newtheorem{itheorem}{Theorem}{\Alph{theorem}}

\newtheorem{iproposition}[itheorem]{Proposition}

\theoremstyle{definition}

\newtheorem{subexample}[subtheorem]{Example}

\newtheorem{example}[theorem]{Example}
\newtheorem{examples}[theorem]{Examples}

\newtheorem{chunk}[theorem]{}
\newtheorem{subchunk}[subtheorem]{}

\theoremstyle{remark}

\newenvironment{bfchunk}{\begin{chunk}\textit}{\end{chunk}}

\newtheorem*{Addendum}{Addendum}

\newtheorem{remark}[theorem]{Remark}

\numberwithin{equation}{theorem}

\newcommand{\add}{\operatorname{add}}
\newcommand{\sadd}[1]{{\mathsf{add}^{\scriptscriptstyle\mathsf
{\Sigma}}(#1)}}

\newcommand{\smd}[1]{{\mathsf{smd}\left(#1\right)}}

\newcommand{\level}[3]{\operatorname{level}_{#1}^{#2}(#3)}
\newcommand{\thickn}[3]{{\mathsf{thick}}_{#2}^{#1}(#3)}
\newcommand{\thick}[2]{{\mathsf{thick}}_{#1}(#2)}
\newcommand{\slength}[3]{{\operatorname{\ell}_{#1}^{#2}(#3)}}
\newcommand{\slol}[3]{{\operatorname{\ell\hspace{-1pt}\ell}_{#1}^{#2}(#3)}}

\newcommand{\prclass}[2]{\operatorname{proj\,class}_{#1}#2}

\newcommand{\cxy}[2][R]{{\operatorname{cx}_{#1}{#2}}}
\newcommand{\filt}[2]{{{#1}^{#2}}}
\newcommand{\Filt}[2]{{({#1}^{#2})}}

\newcommand{\dfm}[1]{{{#1}{}^{\flat}}}

\begin{document}
\title[Homology of perfect complexes]
{Homology of perfect complexes}

\author[L.~L.~Avramov]{Luchezar L.~Avramov}
\address{Department of Mathematics, University of Nebraska, Lincoln, NE 68588, U.S.A.}
\email{avramov@math.unl.edu}
\author[R.-O.~Buchweitz]{Ragnar-Olaf Buchweitz} 
\address{Department of Computer and Mathematical Sciences,\newline\indent
University of Toronto Scarborough, Toronto, ON M1A 1C4, Canada} 
\email{ragnar@utsc.utoronto.ca}
\author[S.~B.~Iyengar]{\\Srikanth B.~Iyengar} 
\address{Department of Mathematics, University of Nebraska, Lincoln, NE 68588, U.S.A.}  
\email{iyengar@math.unl.edu}
\author[C.~Miller]{Claudia Miller} 
\address{Mathematics Department, Syracuse University, Syracuse, NY 13244, U.S.A.}  
\email{clamille@syr.edu}

\thanks {Research partly supported by NSF grants DMS 0201904 and DMS 0803082 (L.L.A); NSERC grant 3-642-114-80 (R.-O.B); NSF grants DMS 0602498 and DMS 0903493 (S.B.I.); NSF grant DMS 0434528 and NSA grant H98230-06-1-0035 (C.M.)}

  \begin{abstract}
It is proved that the sum of the Loewy lengths of the homology modules
of a finite free complex $F$ over a local ring $R$ is bounded below
by a number depending only on  $R$. This result uncovers,  in the
structure of modules of finite projective dimension, obstructions to
realizing $R$ as a closed fiber of some flat local homomorphism. Other
applications include, as special cases, uniform proofs of known results
on free actions  of elementary abelian groups and of tori on finite CW
complexes. The arguments use numerical invariants of objects in general
triangulated categories, introduced here and called levels.  They allow
one  to track,  through changes of triangulated categories, homological
invariants like projective dimension, as well as structural invariants
like Loewy length.  An intermediate result sharpens, with a new proof,
the  New Intersection Theorem for commutative  algebras over  fields.
Under additional hypotheses on the ring $R$ stronger estimates are proved
for Loewy lengths of modules of finite projective dimension.
  \end{abstract}

%\subjclass[2010]{13D02, 13D22, 13D09; 13H10, 18E30, 19J35}

\date{16th October 2009}

\maketitle

\tableofcontents

  \newpage

\section*{Introduction}

We study homological properties of finite free complexes over noetherian
local rings. The results uncover novel links between the structure of
the homology modules of such complexes and conormal modules of the ring.
Their statements incorporate intuition coming from research in algebraic
topology.  The arguments use techniques from commutative algebra,
differential graded homological algebra, and triangulated categories,
and develop new tools for these fields.

Let $(R,\fm,k)$ be a local ring with maximal ideal $\fm$ and residue field
$k$.  One theme that runs through the paper is that when $R$ has
`small' modules of finite projective dimension it is `not far' from 
being regular.  Here we measure the `size' of an $R$-module 
$M$ in terms of its \emph{Loewy length}, defined to be the number
\[
\lol RM = \inf\{i\ge 0\mid \fm^i M= 0\}\,.
\]
When $M$ is non-zero and of finite projective dimension, Loewy length one or two occurs if and only if $R$ is regular or a quadratic hypersurface, respectively.

We provide uniform lower bounds on Loewy lengths of modules of
finite projective dimension in terms of invariants depending only on the
ring $R$.  The first one involves the Castelnuovo-Mumford regularity of
$R$; see Section~\ref{Loewy length of modules of finite projective
dimension} for a definition.

\begin{itheorem}
  \label{ireg}
If $R$ is Gorenstein and its associated graded ring $\gr_{\fm}(R)$ is
Cohen-Macaulay, then each non-zero $R$-module $M$ of finite
projective dimension satisfies
  \[
\lol RM\ge\reg{}{R}+1\,.
  \]
\end{itheorem}

This result applies, for example, when $\gr_{\fm}(R)$ is a graded complete
intersection, that is, $\gr_{\fm}(R) \cong k[x_1,\dots,x_e]/I$ where
$I$ is generated by a homogeneous regular sequence $g_1,\dots,g_c$;
in this case $\reg{}R=\sum_{j=1}^c (\deg g_j-1)$.  Even this special
case of the theorem above was known only for $c=1$, where it was proved
by Ding~\cite{Di:cm}.

Applications of Theorem~\ref{ireg} are restricted by the fact that its
hypothesis is on the associated graded ring, rather than on the ring
of interest.  For instance, $\gr_{\fm}(R)$ need not be {Cohen-Macaulay}
even when $R$ is a local complete intersection.

In the remaining results of this work the hypotheses bear on the structure
of $R$ itself.  Recall that its \emph{embedding dimension} is the number
$\edim R=\rank_k (\fm/\fm^2)$.

\begin{itheorem}
  \label{intro:fibres}
Let $(P,\fp)\to (Q,\fq)$ be a flat local homomorphism and set
$R=Q/\fp Q$.

Every non-zero $R$-module $M$ of finite projective dimension then
satisfies
\[
\lol RM\ge\edim P-\edim Q+\edim R+1\,.
\]
In particular, one has $\fq^{l}\not\subseteq \fp Q$ for 
$l=\edim P-\edim Q+\edim R$.
 \end{itheorem}  

The special case $\edim R=\edim Q$ captures an important aspect of the
theorem: If $\fm^hM=0$ for some $R$-module $M$ of finite projective
dimension, then $R$ admits no embedded (that is to say, with $\fp Q\subseteq \fq^{2}$) 
flat deformation over a base of embedding dimension greater than $h$.

Theorem~\ref{ireg} is proved in Section~\ref{Loewy length of modules
of finite projective dimension}, using invariants of Gorenstein local
rings defined by M.~Auslander and studied by S.~Ding.  The proof of
Theorem~\ref{intro:fibres} is an altogether different affair.  It is
derived from a more general result, which gives information on the
structure of the homology of \emph{finite free complexes}; that is,
complexes of finitely generated free $R$-modules of the form
  \[
F = \quad 0\lra F_s\lra F_{s-1} \lra \cdots \lra F_{t+1}\lra F_t\lra 0
 \]

In our main result we estimate the sum of the Loewy lengths of the
homology modules of $F$ in terms of another invariant of $R$, which
we call the \emph{conormal free rank} and denote $\cnr R$.  When $R$ is
complete equicharacteristic or is essentially of finite type over a field,
$\cnr R$ equals the maximal rank of a free summand of its conormal module;
see Section~\ref{The conormal rank of a local ring}.  The inequality in
the next theorem is sharp, in the sense that equality holds in some cases.

  \begin{itheorem}
\label{intro:loewy}
Every finite free complex $F$ with $\hh F\ne 0$ satisfies an inequality:
\[
\sum_{n\in\BZ}\lol R{\HH nF} \geq \cnr R +1\,.
\]
  \end{itheorem}

Evidently, to apply this result one needs lower bounds for the conormal
free rank.

In the situation of Theorem~\ref{intro:fibres} the free $R$-module $\fp
Q/(\fp Q)^{2}$ is the conormal module of the surjection $Q\to R$, which
can be used to show that $\cnr R$ is not less than $l$.  This allows
one to deduce Theorem~\ref{intro:fibres} from Theorem~\ref{intro:loewy}.

One case when the value of $\cnr R$ is known is when $R$ is complete
intersection: it is equal to the codimension of $R$.  Remarkably, over
such rings one can bound Loewy length of homology for every homologically
finite complex, see Theorem~\ref{intro:ciall}.

The proof of Theorem~\ref{intro:loewy} is presented in Section~\ref{Loewy
length of homology of perfect complexes}.  It draws on the results
of sections~\ref{Levels in triangulated categories} through~\ref{DG
algebra models for Koszul complexes} and involves a number of concepts
and techniques that are not traditional for commutative algebra.

The first move is to replace complexes over rings with DG (that is,
differential graded) modules over DG algebras.  Such a procedure
was developed by Avramov to study homological invariants of rings and
modules; see \cite{Av:ifr} for a survey.  It utilizes the possibility of
adjusting the algebraic properties of a DG algebra by replacing it with
a quasi-isomorphic one, while replacing its derived category with an
equivalent category.  However, the situation here is different because
structural properties, such as Loewy length, depend on the underlying
graded algebra of a DG algebra and need not be preserved by equivalences
of derived categories.

A crucial new idea is to bound Loewy length by numerical invariants of
objects of the derived category of the ring that behave predictably under
applications of exact functors.  Their introduction is motivated in part
by work of Dwyer, Greenlees, and Iyengar \cite{DGI2}.  These authors
transported from homotopy theory into commutative algebra the concept
of \emph{building} an object $X$ in a triangulated category $\ct$
from some fixed object $C$ of $\ct$.  Here we define a number, $\level
{\ct}{C}{X}$, that we call the \emph{level} of $X$ with respect to $C$
in $\ct$.  It measures the number of \emph{extensions} needed in the
`building process'.  Further suggestions that such a notion might be
useful came from the papers of D.~Christensen \cite{Ch}, Bondal and
Van den Bergh \cite{BV}, and Rouquier~\cite {Rq:rd,Rq:dim} dealing with
dimensions of triangulated categories.

To show how these ideas fit together we sketch an outline of the proof
of Theorem~\ref{intro:loewy}.  It also serves as an overview of the
content of the paper.

In Section~\ref{Levels in triangulated categories} we define levels and
record their elementary properties.

In Section~\ref{Levels of DG modules} we specialize to the case of
the derived category $\dcat A$ of DG modules over a DG algebra $A$;
for simplicity, we write $\level {A}{C}{X}$ in place of $\level{\dcat
A}{C}{X}$.  Two levels over $A$ play a special role in this work.

Levels with respect to the DG $A$-module $A$ extend the concept of
projective dimension  from modules over a ring to DG modules over a DG
algebra.  For instance, for every complex $F$ of finite free $R$-modules
the definitions give an inequality
  \begin{equation}
   \label{ipd}
\level RRF\leq\card\{n\in\BN\mid F_n\ne0\}\,.
  \tag{$*$}
   \end{equation}

A structure theorem for DG $A$-modules of finite $A$-level is proved in
Section~\ref{Perfect DG modules}.  In Section~\ref{A New Intersection
Theorem for DG  algebras} it is used, along with the main theorem of
\cite{dm}, to prove the result below.  In view of \eqref{ipd}, the
inequality on the right generalizes and sharpens the classical New
Intersection Theorem for commutative noetherian algebras over fields.

  \begin{itheorem}
 \label{intro:plevel}
Let $A$ be a DG algebra with zero differential, let $M$ be a DG module
over $A$, and let $I$ denote the annihilator of $\hh M = \bigoplus_{n\in\BZ}\HH
nM$ in the ring $A^\flat=\bigoplus_{n\in\BZ} A_n$.

When $A^\flat$ is a commutative noetherian algebra over a field one has
 \[
\pd_A\hh M+1 \ge \level AAM\geq \height I+1  \,.
 \]
  \end{itheorem}

Of major importance here is also the level with respect to a semi-simple
DG $A$-module $k$; its behavior is akin to that of Loewy length of
modules over rings.  The basic properties of $\level Ak-$ are derived
in Section~\ref{Levels and semi-simplicity}.  In particular, we prove

  \begin{iproposition}
\label{intro:length}
If $M$ is a complex of $R$-modules and $\hh M$ has finite length, then
\[
\sum_{n\in\BZ}\lol R{\HH nM} \geq \level RkM
\ge\max_{n\in\BZ}\{\lol R{\HH nM}\}\,.
\]
  \end{iproposition}

The next stage in the proof of Theorem~\ref{intro:loewy} is to construct
a chain of exact functors
   \begin{equation*}
\xymatrixcolsep{1.5pc}
    \xymatrix{
\dcat R \ar@{->}[r]^{\ST}
     &\dcat K\ar@{->}[r]^-{\SJ}_-\equiv
       &\dcat{\Lambda\otimes_kB}\ar@{->}[r]^-{\SII}
        &\dcat \Lambda
}
  \end{equation*}
of derived categories of DG modules, where $\equiv$ flags an equivalence.
The DG algebras are:  the ring $R$ viewed as a DG algebra concentrated
in degree $0$; the Koszul complex $K$ on a minimal generating set for
$\fm$; a DG $k$-algebra $B$ with $\rank_kB<\infty$; and an exterior
algebra $\Lambda=k\langle\xi_1,\dots,\xi_c\rangle$ with $\deg \xi_i=1$,
zero differential, $c=\cnr R$.

Two functors are easy to describe: $\ST(-)=K\otimes_R-$ and
$\SII$ is the forgetful functor defined by the  canonical morphism
$\Lambda\to\Lambda\otimes_kB$.  On the other hand, the construction of the
equivalence of categories $\SJ$ takes up all of Section~\ref{DG algebra
models for Koszul complexes}, where we also prove that the DG module
$\SJ\ST(k)\in\dcat{\Lambda\otimes_kB}$ is a direct sum of suspensions
of copies of $k$; as a consequence, levels with respect to $\SJ\ST(k)$
are equal to $k$-levels.  This fact, together with the formal property
that functors do not raise levels and equivalences preserve them,
justify all but the initial step in the following sequence:
  \begin{align*}
\sum_{n\in\BZ}\lol R{\HH nF}
   &\geq \level{R}{k}{F} \geq\level{K}{\ST(k)}{\ST(F)}\\
   &=\level{\Lambda\otimes_kB}{\SJ\ST(k)}{\SJ\ST(F)}
   =\level{\Lambda\otimes_kB}{k}{\SJ\ST(F)} \geq\level{\Lambda}{k}{N}\,,
  \end{align*}
where we have set $N=\SII\SJ\ST(F)$.  The first inequality is provided
by Proposition~\ref{intro:length}.

To finish the proof of Theorem~\ref{intro:loewy} one needs an estimate
for $\level{\Lambda}{k}{N}$.  The key to obtaining one is to show that
$N$ is isomorphic in $\dcat{\Lambda}$ to a DG module with finite free
underlying graded $\Lambda$-module.  To do this we first remark that
formula \eqref{ipd} implies that the finite free complex $F$ has finite
$R$-level, then use results on persistence of levels to show that $N$
has finite $\Lambda$-level, and finally prove that the last condition
is equivalent to $N$ being finite free as a graded $\Lambda$-module.

It remains to use the equality in the following result of independent
interest:

  \begin{itheorem}
\label{intro:bgg}
If $N$ is a DG $\Lambda$-module with finite free underlying
graded $\Lambda$-module and with $\hh N\ne 0$, then one has
\[
\card\{n\mid \HH nN\ne 0\}\ge \level{\Lambda}kN=c+1\,.
\]
  \end{itheorem}

The proof of Theorem \ref{intro:bgg} itself consists of two independent
steps.  The first one is to reduce the computation of $\level{\Lambda}kN$
to that of $\level SSM$, where $S$ is a polynomial ring in $c$
indeterminates over $k$, and $M$ is a DG $S$-module with $\rank_k\hh M$
finite.  For this we use a DG version of the Bernstein-Gelfand-Gelfand
equivalence, presented in Section~\ref{Perfect DG modules over exterior
algebras}.  Theorem \ref{intro:plevel} then gives $\level SSM=c+1$.

All the threads of the delicate proof of Theorem~\ref{intro:loewy} come
together in Section~\ref{Loewy length of   homology of perfect complexes}.
In the final count, the argument hinges on the ability to simultaneously
track numerical invariants, both structural and homological, under the
actions of various functors.  Its success attests to the remarkable
versatility of the concept of level.

Theorem~\ref{intro:loewy} is restricted to finite free complexes.
Using results from \cite{Av:vpd}, over certain rings we extend it to a
statement about all complexes with finite homology:

\begin{itheorem}
\label{intro:ciall}
If $R$ is a complete intersection local ring and $M$ is a complex of 
$R$-modules with $\hh M$ finite and non-zero, then one has an inequality
\[
\sum_{n\in\BZ}\lol R{\HH nM}  \geq \codim R - \cxy M + 1\,.
\]
\end{itheorem}

This is proved in Section~\ref{Complete intersection local rings}. The
number $\cxy M$, known as the \emph{complexity} of $M$, is the least
non-negative integer $d$ such that the ranks of the free modules in a
minimal free resolution of $M$ are bounded by a polynomial of degree $d-1$.
Theorem~\ref{intro:ciall} can also be deduced from a strengthening of
Theorem~\ref{intro:bgg} that covers all DG $\Lambda$-modules with finite
underlying graded module; see \cite{AI:leeds}.

The last two theorems above have antecedents in the study of the homology
of a finite CW complex $X$ with an action of an elementary abelian group
or of a torus.

When $G$ is an elementary abelian $p$-group and $X$ has a $G$-equivariant cellular decomposition, its cellular
chain complex with coefficients in $\BF_{p}$ is finite free over the ring $T\cong \BF_p[x_1,\dots,x_c]/(x_1^p,\dots,x_c^p)$. Carlsson \cite{Ca:inv} established that
  \[
\sum_ n\lol T\HH n{X\,;\BF_p}\ge c+1
  \]
holds by proving Theorem~\ref{intro:ciall} for $M$ a finite free complex over
$T$ and $p=2$.  Allday, Baumgartner, and Puppe, see \cite{AP}, proved
Theorem~\ref{intro:ciall} for odd $p$ and for all complexes of $T$-modules
with finite homology; the arguments depend on the parity of $p$.

When the induced action of $G$ on $\HH *{X\,;\BF_p}$ is trivial the
inequality above simply states that $X$ has at least $c+1$ non-trivial
homology groups with coefficients in $\BF_p$.  Allday and Puppe
\cite{AP} proved a similar estimate for almost free torus actions:
  \[
\card\{n\mid \HH n{X\,;\BQ}\ne 0\}\ge c+1\,.
  \]
The algebraic core of their proof is a property of DG modules over a 
polynomial ring $S$ in $c$ indeterminates over $\BQ$, which is implied 
by Theorem~\ref{intro:plevel}.

The original proofs of the theorems on group actions heavily depend on
the structure of the rings  $S$ and $T$.  Our results demonstrate that 
these theorems are manifestations of general properties of complexes 
over commutative noetherian rings.

\section{Loewy length of modules of finite projective dimension}
\label{Loewy length of modules of finite projective dimension}

Let $R$ be a commutative local noetherian ring with maximal ideal $\fm$
and residue field $k=R/\fm$.  The \emph{Loewy length} of an $R$-module
$M$ is the number
   \[
\lol RM=\inf\{n\in\BN \mid \fm^nM=0\}\,.
   \]
   When $M$ is finitely generated, $\lol RM$ is finite if and only if its
   (Jordan-H\"older) \emph{length}, denoted $\length_RM$, is finite.  Often the Loewy length of $M$
   carries more structural information than does its length.

Let $G$ denote the associated graded ring $\gr_{\fm}(R)$; thus, $G_{0}=k$ and $G$ is generated over $G_{0}$ by $G_{1}$.  Choose a presentation $G\cong P/J$ with $P=k[x_1,\dots,x_e]$, a polynomial ring in indeterminates $x_i$ of degree $1$, and $J$ a homogeneous ideal.  Set
    \[
\reg{}R=\sup\{j\mid \Tor iPG{P/P_{\ges 1}}_{i+j}\ne0\}\,.
    \]
In other words, $\reg{}R$ is the \emph{Castelnuovo-Mumford regularity},
$\reg PG$, of the graded $P$-module $G$; it is independent of the choice
of $P$, see e.g.~\c{S}ega~\cite[1.4]{Sg}.

The \emph{order} of the ring $R$ is given for singular $R$ by the formula
  \[
\ord R=\inf\left\{n\in\BN\,\left\vert\,\length_R(R/\fm^{n+1})
<\binom{n+e}{e}\right\}\right.
  \]
where $e=\edim R$; when the ring $R$ is regular we set $\ord R=1$.

The next result contains Theorem \ref{ireg} from the introduction. It
applies, in particular, when $R$ is the localization or the completion
of a standard graded Gorenstein ring $G$ at the maximal ideal
$\bigoplus_{n\ges1}G_n$.

  \begin{theorem}
  \label{bound:modules}
If $R$ is Gorenstein and the associated graded ring $\gr_{\fm}(R)$ is
  Cohen-Macaulay, then for each non-zero $R$-module $M$ of finite projective dimension
  \[
\lol RM \geq\reg{}R+1\ge\ord R\,.
   \]
When $k$ is infinite the first inequality becomes an equality for
some such module $M$.
  \end{theorem}

The proof uses numerical invariants introduced by M.~Auslander.

Let $R$ be a Gorenstein local ring and $M$ a finite $R$-module. Let
$M^{\operatorname{cm}}$ denote the sum of all submodules
$\lambda(L)\subseteq M$, when $L$ ranges over all maximal Cohen-Macaulay
$R$-modules with no non-zero free direct summand, and $\lambda$ ranges
over all $R$-linear homomorphisms $L\to M$.  The minimal number of
generators of the $R$-module $M/M^{\operatorname{cm}}$ is called the
\emph{delta invariant} of $M$ and is denoted $\delta_R(M)$.  

The basic properties of $\delta_R(-)$, due to Auslander, are collected
in the next result. Their proofs in the literature are scattered and
some use alternative characterizations, so we provide complete details.

  \begin{lemma}
  \label{index}
Let $M$, $N$ be finite modules over a Gorenstein local ring $(R,\fm,k)$.

The following (in)equalities then hold:
    \begin{align}
 \label{xi:additive}
&\delta_R(M)=\delta_R(M')+\delta_R(M'')
\quad\text{when $M=M'\oplus M''$}\,.
\\
 \label{xi:monotonic}
&\delta_R(M)\ge\delta_R(N)
\quad\text{when $M \to N$ is a surjective homomorphism}.
\\
 \label{xi:maximal}
&\delta_R(M)=\rank_k(M/\fm M)\quad\text{when $\pd_RM<\infty$}\,.
\\
 \label{xi:index}
&\delta_R(R/\fm^n)=1\quad\text{for all $n\gg0$}\,.
\\
 \label{xi:regular}
&\delta_R(k)=1\quad\text{if and only if $R$ is regular}\,.
  \end{align}
    \end{lemma}

  \begin{proof}
The first two assertions are evident. 

For the third one, it suffices to prove that if $L$ is a maximal
Cohen-Macaulay module and $\lambda\col L\to M$ is a homomorphism of
$R$-modules with $\lambda(L)\nsubseteq\fm M$, then $L$ has a non-zero
free summand.  Consider an exact sequence
\[
0\lra E\lra F\xra{\ \pi\ } M\lra 0
\]
with $F$ free and $E\subseteq \fm F$.  By \cite[3.3.10(d)]{BH} and
induction on $\pd_RE$, one gets $\Ext 1RLE=0$, so the following
map is surjective:
\[
\Hom RL\pi\col \Hom RLF\lra \Hom RLM\,.
\]
Thus, there is a homomorphism $\rho\col L\to F$ with $\pi\rho=\lambda$.
Choosing $x\in L$ with $\lambda(x)$ not in $\fm M$ we get $y=
\rho(x)\notin\fm F$, so $Ry$ is a non-zero free direct summand of $F$. The
composition $L\to F\to Ry$ is then surjective, as desired.

To prove the fourth formula, choose a maximal $R$-regular sequence
$\bsx$. Since $\length_R(R/R\bsx)$ is finite, for all $n\gg0$ one
has $\fm^n\subseteq R\bsx$. Since $\pd_R(R/R\bsx)$ is finite the
surjection $R/\fm^n\to R/R\bsx$ implies $\delta_R(R/\fm^n)\ge1$, see
\eqref{xi:maximal}. It remains to note that any finite $R$-module $M$ 
satisfies $\delta_{R}(M)\leq \rank_{k}(M/\fm M)$.

As to the last one, $R$ regular implies $\pd_Rk$ finite, hence
$\delta_R(k)=1$ by \eqref{xi:maximal}. If $R$ is not regular, then
it has a maximal Cohen-Macaulay $R$-module $L$ with no non-zero free
summands. Picking any surjection $L\to k$ one gets $\delta_R(k)=0$.
  \end{proof}

The \emph{index} of $R$ is defined by Auslander to be the number
\[
  \indx R=\inf\left\{n\in\BN \mid \delta_R(R/\fm^n)\ge1\right\}\,.
\]
It is a positive integer, by \eqref{xi:index}, and equals $1$ if and
only if $R$ is regular, by \eqref{xi:regular}.  

The next result is the first step in the proof of
Theorem~\ref{bound:modules}.

  \begin{lemma}
 \label{index:bound}
Let $R$ be a Gorenstein local ring and $M$ a finite, non-zero $R$-module.
If the projective dimension of $M$ is  finite, then one has
  \[
\lol RM \geq \indx R\,.
  \]
    \end{lemma}

  \begin{proof}
We may assume $\lol RM=l<\infty$.  For $r=\rank_k(M/\fm M)$ there
is then a surjection $(R/\fm^l){}^r \to M$.  The formulas in
Lemma \ref{index} yield (in)equalities
  \[
r\cdot\delta_R(R/\fm^l)
=\delta_R\left( (R/\fm^l){}^r\right)\geq\delta_R(M)=r\,.
  \]
They imply $\delta_R(R/\fm^l)\ge1$, which means $l\ge\indx R$.
    \end{proof}

To get lower bounds for $\indx R$ we use a result of Ding.  In order
to state it we recall that a sequence $x_1,\dots,x_d$ of elements of
$\fm$ is said to be \emph{super-regular} if their initial forms in the
associated graded ring $\gr_{\fm}(R)$ form a regular sequence.

  \begin{chunk} \label{index:ding}
When $R$ is Gorenstein, $\operatorname{gr}_{\fm}(R)$ is
Cohen-Macaulay, and $\bsx$ is a super-regular sequence
of length $\dim R$ in $\fm\smallsetminus\fm^2$, the \emph{proof} of
\cite[2.1]{Di:cm} yields an equality:
    \[
\indx R=\lol R{(R/R\bsx)}\,.
    \]
  \end{chunk}

  \begin{lemma}
    \label{index:regularity}
If the ring $R$ is Gorenstein and $\gr_{\fm}(R)$ is Cohen-Macaulay, then
\[
\indx R = \reg{}R+1\,.
  \]
    \end{lemma}

  \begin{proof}
Set $G=\gr_{\fm} (R)$, and choose a presentation $G\cong P/J$, where $P$ is a polynomial ring over $k$, generated by finitely many indeterminates of degree $1$.

After a flat base change one may assume that $k$ is infinite.  One can
then find in $P_1$ a sequence $\bsy=y_1,\dots,y_d$, with $d=\dim G$,
which is both $G$-regular and $P$-regular.  For $i=1,\dots,d$ choose
$x_i\in\fm\smallsetminus\fm^2$ with initial form $y_i\in G_1$.  The
sequence $\bsx=x_1,\dots,x_d$ is then super-regular, and one has a
chain of equalities
    \begin{align*}
   \indx R &=\lol R{(R/R\bsx)}
\\
    &=\inf\{n\in\BN \mid (G_1)^n\big(\gr_\fm(R/R\bsx)\big)=0\}
\\
    &=\inf\{n\in\BN \mid (G_1)^n(G/G\bsy)=0\}
\\
    &=\reg {P}{(G/G\bsy)}+1
\\
    &=\reg{P}{G}+1
     \end{align*}
that come from \ref{index:ding}, the definition of Loewy length,
the isomorphism $\gr_\fm(R/R\bsx)\cong G/G\bsy$, see \cite[0.1]{Sa},
and standard properties of regularity.
\end{proof}

 One last observation is needed before proving the theorem.

  \begin{chunk} 
  \label{index:order}
For every local ring $R$ one has $\reg{}R\ge\ord R-1$.

Indeed, let $\wh R\cong Q/I$ be a minimal regular presentation; see
\ref{cohen}.  By passing to associated graded rings one gets from it an
isomorphism $G\cong P/J$, where $P$  is a polynomial ring in $e=\edim
R$ indeterminates.  One has $\rank_k P/(P_1)^{n+1}=\binom{n+e}e$, so
$\ord R$ is equal to the least degree of a non-zero element in $J$.
The isomorphism $\Tor 1PGk\cong J\otimes_Pk$ of graded vector spaces
then yields the desired inequality.
    \end{chunk}

  \begin{proof}[Proof of Theorem \emph{\ref{bound:modules}}]
Lemma~\ref{index:bound}, Lemma~\ref{index:regularity}, and
\ref{index:order}, imply
  \[
\lol RM\ge\reg {}R+1\ge\ord R\,.
  \]
When $k$ is infinite the $R$-module $M'=R/R\bsx$ from the proof
of Lemma~\ref{index:regularity} has $\pd_RM'<\infty$.  The computation
there yields $\lol R{M'}=\reg PG+1$.
  \end{proof}

In many situations good estimates of the regularity of graded rings are
known.  The next result provides one that seems to be new. It involves
conormal modules, which also appear in our main result, see 
Theorem~\ref{conormal:main}.

  \begin{proposition}
   \label{reg:conormal}
Suppose $G\cong P/J$ with $P$ a standard graded polynomial ring
over $k$ and $J$ a homogeneous ideal.  

If $J/J^2$ is isomorphic to $\big(\bigoplus_{j=1}^c G(-n_j)\big) \oplus C$ 
as graded $G$-modules, then 
  \[
\reg {P}G \geq
  \begin{cases}
\sum_{j=1}^c (n_j-1)&\text{when $C=0$}\,,\\ \sum_{j=1}^c
(n_j-1)+(m-1)&\text{otherwise}
  \end{cases}
\] 
where $m$ is the least degree of a non-zero element in $J$.
\end{proposition}

  \begin{proof}
Let $K$ be the Koszul complex on a $k$-basis of $G_1$, considered
as a complex of graded $G$-modules.  For each $i$ one then has $\Tor
iPGk\cong\HH iK$ as graded $k$-vector spaces.  The arguments in the
proofs of \cite[(2.3) and (2.1)]{Iy}, see Theorem~\ref{koszul:structure}, 
show that $K$ is quasi-isomorphic to a complex of graded $k$-vector spaces
$V\otimes_kW$, where
  \[
V=\quad\cdots
\xra{\ 0\ } {\textstyle\bigwedge}^{i+1}\left(\bigoplus_{j=1}^c k(-n_j) 
\right)
\xra{\ 0\ } {\textstyle\bigwedge}^i\left(\bigoplus_{j=1}^c k(-n_j) 
\right)\xra{\ 0\ }\cdots
  \]
and $W$ satisfies $W_0=k$ and $\dd(W_1)=0=\dd(W_2)$.  The isomorphisms
  \[
V_c\cong k(-n)
\quad\text{and}\quad
\hh{V\otimes_kW}\cong V\otimes_k\hh W\,,
  \]
where $n=\sum_{j=1}^cn_j$, define injective linear maps of
graded vector spaces
  \begin{align}
  \label{mono1}
k(-n)\otimes_kW_0&\lra\Tor cPGk
\\
  \label{mono2}
k(-n)\otimes_kW_1&\lra\Tor{c+1}PGk
  \end{align}

{F}rom \eqref{mono1} one gets $\reg PG\ge n-c$, which is the first
inequality of the proposition.  Set $q=\inf\{j\in\BZ\var(W_1)_j\ne0\}$.
The isomorphisms
  \[
\left(\bigoplus_{j=1}^c k(-n_j)\right)\oplus(C\otimes_Pk)\cong
J\otimes_Pk\cong
\Tor1PGk\cong
V_1\oplus W_1
  \]
of graded $k$-vector spaces yield $W_1\cong C\otimes_Pk\subseteq J 
\otimes_Pk$.
Thus, $C\ne0$ implies $q\ge m$, so \eqref{mono2} gives
$\reg {}G\ge(m+n)-(c+1)$.  This is the second desired inequality.
  \end{proof}

The theorem is readily applicable to a class of complete intersection
rings that includes all hypersurface rings.   The  notation is that
of Theorem~\ref{bound:modules}.

  \begin{example}
    \label{index:sci}
Suppose the local ring $(R,\fm,k)$ is a \emph{strict complete intersection} of type$(n_1,\dots,n_c)$, meaning that for some isomorphism $\gr_{\fm}(R)\cong P/J$ the ideal $J$ can be generated by a homogeneous $P$-regular sequence $\bsg=g_1,\dots,g_c$ with $\deg g_j=n_{j}$.

For every $R$-module $M$ of finite projective dimension one has
\[
     \lol RM\ge \sum_{j=1}^c (n_j - 1) + 1\,.
\]

Indeed, the hypothesis implies an isomorphism $J/J^2\cong \bigoplus_{j=1}^c\gr_{\fm}(R)(-n_j)$ of graded modules over $\gr_{\fm}(R)$, so Theorem~\ref{bound:modules} and Proposition~\ref{reg:conormal} apply.
    \end{example}

\section{Levels in triangulated categories}
\label{Levels in triangulated categories}

To any object, say $X$, in a triangulated category we associate a
numerical invariant that measures the minimal number of steps necessary
to `build' $X$ out of objects belonging to some fixed reference class.
Its definition uses techniques developed by Beilinson, Bernstein,
and Deligne \cite[\S1.3]{BBD}, Bondal and Van den Bergh \cite[\S1]{BV},
D. Christensen \cite{Ch}, and Rouquier \cite{Rq:dim}, and  is partly
motivated by the work of these authors.  Where our approach differs from
theirs is in the focus on properties of individual objects, rather than
on global invariants of triangulated categories.

In this section we review the constructions allowed in the building
process alluded to above, define levels, and record some general features.
Only basic properties of triangulated categories are needed, and they
can be found in Krause's succinct exposition \cite[\S\S 1--3]{Kr};
for more details we refer the reader to Neeman's book \cite{Ne}.

We use $\shift$ to denote the suspension functor in a triangulated category.
Let $\cs$ be a subcategory---always assumed non-empty---of a triangulated
category $\ct$.  We say that $\cs$ is \emph{strict} if it is closed under
isomorphisms in $\ct$; it is \emph{full} if every morphism in $\ct$
between objects in $\cs$ is contained in $\cs$. A strict full subcategory $\cs$ is \emph{thick} if
it is additive, closed under direct summands, and in any exact triangle
$L\to M\to N\to\shift L$ in $\ct$ when two of the objects $L,M,N$ are in $\cs$
so is the third. Note that every thick subcategory is  triangulated.

  \begin{bfchunk}{Operations on subcategories.}
Let $\ct$ be a triangulated category and $\ca$ a subcategory of $\ct$.
We define several closure operations on $\ca$ in $\ct$.

  \begin{subchunk}
\label{Additive closures}
The intersection of all strict and full subcategories of $\ct$
that contain $\ca$ and are closed under finite direct sums and all
suspensions is denoted $ \sadd{\ca}$; in \cite{BV} this subcategory is
denoted $\add{\ca}$.  
  \end{subchunk}

  \begin{subchunk}
\label{Karoubization}
The intersection of all full subcategories of $\ct$ that contain $\ca $ and are closed under retracts (equivalently, isomorphisms and direct summands) is denoted $\smd\ca$.
    \end{subchunk}

  \begin{subchunk}
\label{star}
Given strict and full triangulated subcategories $\ca$ and $\cb$ of
$ \ct$, let $\ca\star\cb$ be the full subcategory whose objects are
described as follows:
   \[
\ca\star\cb=\left\{M\in \ct\left | 
  \begin{gathered}
     \text{there is an exact triangle}
\\
      L\to M\to N\to\shift L
\\
     \text{with $L\in\ca$ and $N\in\cb$}
  \end{gathered} \right\}\right.
    \]
This subcategory is strict.  For every strict and full subcategory  $
\cc$ of $\ct$ one has
   \[
    \ca\star(\cb\star \cc) = (\ca\star \cb)\star \cc
   \]
see \cite[(1.3.10)]{BBD}.  Thus, the following notation is unambiguous:
   \[
   \ca^{\star n} =
   \begin{cases}
    \{0\}&\text{ for } n=0\,;
\\
   \ca&\text{ for } n=1\,;
\\
    \overset{n\text{ copies}}{\overbrace{\ca\star\cdots\star\ca}}
    &\text{ for } n\geq 2\,,
   \end{cases}
\]
We refer to the objects of $\ca^{\star n}$ as \emph{$(n-1)$-fold
extensions} of objects from $\ca$.  
  \end{subchunk}

  \begin{subchunk}
\label{thick:definition}
Let $\cc$ be a subcategory of $\ct$.  The intersection
of all thick subcategories of $\ct$ containing $\cc$ is itself a thick
subcategory, called the \emph{thick closure} of $\cc$ in $\ct$; in this
paper it is denoted by $\thick{\ct}{\cc}$.
    \end{subchunk}
     \end{bfchunk}

Properties of objects in a subcategory often propagate to its thick
closure:

  \begin{subchunk}
\label{levels:inheritance}
Let $\mathcal P$ denote a property of objects in $\ct$, and assume
that the full subcategory consisting of the objects with property
$\mathcal P$ is thick.

If each $C\in \cc$ has property $\mathcal P$, then so does
every object in $\thick{\ct}\cc$.
    \end{subchunk}

One can approximate $\thick{\ct}{\cc}$ `from below' by
a process used in \cite{Ch}, \cite{BV}, \cite{Rq:dim}.

  \begin{bfchunk}{Thickenings.}
\label{thickenings}
Let $\cc$ be a subcategory of $\ct$.

For each $n\in\BN$ the $n$th \emph{thickening} of $\cc$ is the
full subcategory with objects
\begin{equation*}
     \thickn n{\ct}{\cc} =
\begin{cases}
    \{0\} & \text{if }n=0 \\
\smd{\sadd{\cc}}  & \text{if }n=1 \\
\smd{\thickn {n-1}{\ct}{\cc}\star\thickn 1{\ct}{\cc}} &\text{if }n\geq 2
  \end{cases}
    \end{equation*}
This subcategory appears implicitly in \cite{Ch} and explicitly in
\cite{BV}, where it is denoted $\la\cc\ra_{n}$ but is not named.
It is closed under suspensions, finite direct sums, and retracts, but
it is not necessarily closed under formation of exact triangles.

  \begin{subchunk}
\label{layers:symmetry}
The equality below, see \cite[\S 2.1, p.5]{BV}, provides
an alternative description:
\begin{equation*}
\thickn {n}{\ct}{\cc} =\smd{\sadd{\cc}^{\star n}}\,.
  \end{equation*}
In words: The objects in the $n$th thickening of $\cc$ are retracts
of $(n-1)$-fold extensions of objects in $\sadd\cc$.  Thus, the $n$th
thickening can also be built out of the $(n-1)$st one by gluing
objects to the left.
    \end{subchunk}

  \begin{subchunk}
\label{thickenings:relative}
If $\cc$ is contained in some thick subcategory $\cs$ of $\ct$, then
it follows from the definitions that for each integer $n\ge0$
one has $\thickn n{\cs}{\cc}=\thickn n{\ct}{\cc}$.
    \end{subchunk}

  \begin{subchunk}
\label{layers:opposite}
Letting $(-)^\op$ denote passage to the opposite category, from
\ref{layers:symmetry} one gets
    \begin{equation*}
     \thickn n{\ct}{\cc}^{\op} = \thickn n{\ct^\op}{\cc^{\op}}
  \end{equation*}
    \end{subchunk}

  \begin{subchunk}
\label{layers:filtration}
The thickenings of $\cc$ provide a natural filtration of its thick
closure in $\ct$:
\[
\{0\}=\thickn 0{\ct}{\cc}\subseteq \thickn 1{\ct}{\cc}\subseteq\cdots
\subseteq\bigcup_{n\in\BZ} \thickn n{\ct}{\cc} = \thick\ct\cc
\]
  \end{subchunk}
    \end{bfchunk}

In \cite {BV,Rq:dim} the filtration above is used to define the dimension
of $\ct$ as the infimum of the integers $d\ge0$ with the property that
$\thickn {d+1}{\ct}C=\ct$ for some object $C\in\ct$.  Here we derive
from the filtration numerical invariants of the objects of $\ct$.

  \begin{bfchunk}{Levels.}
\label{levels}
Let $\cc$ be a subcategory of  $\ct$.

To each object $M$ in $\ct$ we associate the number
\[
\level{\ct}{\cc}M = \inf\{n\in \BN\mid M\in \thickn n{\ct}{\cc}\}
\]
and call it the \emph{$\cc$-level of $M$}.  It measures the number
of steps required to build $M$ out of $\sadd {\cc}$.  Evidently,
$\level{\ct}{\cc}M<\infty $ is equivalent to $M\in\thick \ct\cc$.

When $\cc$ consists of a single object, $C$, we write $\level{\ct}{C}M $
in place of $\level{\ct}{\{C\}}M$.
  \end{bfchunk}

  \begin{lemma}
    \label{levels:properties}
For each object $M\in\ct$ the following statements hold.
    \begin{enumerate}[\quad\rm(1)]
    \item
$\level{\ct}{\cc}{\shift^iM}=\level{\ct}{\cc}M$ for all $i\in\BZ$.
    \item
If $N\to M\to P\to \shift N$ is an exact triangle in $\ct$, then
    \[
\level{\ct}{\cc}M \leq \level{\ct}{\cc}N + \level{\ct}{\cc}P
    \]
    \item
$\level{\ct}{\cc}{M'\oplus M''} =
\max\{\level{\ct}{\cc}{M'},\level{\ct}{\cc}{M''}\}$
    \item
If $\cc$ is contained in a thick subcategory $\cs$ of $\ct$, then
    \[
\level{\cs}{\cc}M =\level{\ct}{\cc}M\,.
    \]
    \item
$\level{\ct}{\cc}M =\level{\ct^\op}{\cc^{\op}}M$.
    \item
If $\SF\col \ct\to \cu$, respectively, $\SF\col \ct^\op \to \cu$,
is an exact functor of triangulated categories, then the following
inequality holds:
    \[
\level{\cu}{\SF(\cc)}{\SF(M)}  \leq \level{\ct}{\cc}M
    \]
Equality holds if there is a functor $\SG\col \cu\to\ct$ with
$\SG\SF\simeq \idmap^{\ct}$, respectively, $\SG\col \cu\to\ct^{\op}$ with
$\SG\SF\simeq \idmap^{\ct^{\op}}$; in particular,  when $\SF$ is an equivalence.
  \end{enumerate}
    \end{lemma}

  \begin{proof}
The (in)equalities in (1), (2), and (3) come from the definition
of levels.  Parts (4) and (5) come from \ref{thickenings:relative}
and \ref{layers:opposite}, respectively.  Part (6) clearly holds when
$\SF$ is defined on $\ct$; given (5), the case when it is defined on
$\ct^{\op}$ follows.
    \end{proof}

We close the general discussion of levels with some special  features
of these invariants in the case of the derived category of an
associative ring.

  \begin{bfchunk}{Levels of complexes of modules.}
\label{Levels of complexes of modules}
Let $R$ be an associative ring.

In this paper $R$ always acts on its modules from the left.  Our gradings
are usually `homological'.  Thus, complexes of $R$-modules have the form
  \[
M=\quad \cdots \lra M_{n+1} \xra{\dd_{n+1}} M_{n} \xra{\ \dd_n\ }
M_{n-1}\lra \cdots
  \]
For each integer $d$ the $d$th \emph{suspension} of $M$ is the complex
defined by
  \[
(\shift^d M)_n = M_{n-d}
  \quad\text{and}\quad
\dd^{\shift^d M} = (-1)^d\dd^M\,.
  \]
We identify $R$-modules with complexes concentrated in degree
$0$, and graded $R$-modules with complexes with zero differential.
Quasi-isomorphisms are morphisms of complexes that induce isomorphisms
in homology; they are flagged with the symbol $\simeq$, while $\cong$
is reserved for isomorphisms.

A complex of $R$-modules $M$ is \emph{bounded above}, respectively,
\emph{bounded below}, when $M_i=0$ holds for all $i\gg0$, respectively,
for all $i\ll0$; it is \emph{bounded} when both conditions hold; it is
\emph{perfect} if it is quasi-isomorphic to a bounded complex of finite
projective $R$-modules.

Let $\dcat R$ denote the unbounded derived category of the category
of $R$-modules; see \cite[3.2]{Kr} for a description of ${\dcat R}$
as a triangulated category.  We extend the use of $\simeq$ to mark also
isomorphisms in $\dcat R$.

  \begin{subchunk}
\label{dropD}
Let $\cs$ be a thick subcategory of $\dcat R$, such as the derived
category of homologically bounded above (respectively, bounded below,
bounded, or perfect) complexes.  For each subcategory $\cc$ in
$\cs$, each $n\ge0$, and each $M\in\cs$ one has
    \[
\thickn n{\cs}{\cc}=\thickn n{\dcat R}{\cc} \qquad\text{and}\qquad
\level{\cs}{\cc}{M}=\level{\dcat R}{\cc}{M}
    \]
by \ref{thickenings:relative} and Lemma~\ref{levels:properties}(4).
We use the abbreviations
    \[
\thickn n{R}{\cc}=\thickn n{\dcat R}{\cc} \qquad\text{and}\qquad
\level{R}{\cc}{M}=\level{\dcat R}{\cc}{M}
    \]
  \end{subchunk}

Levels of a complex have useful relations to the corresponding levels
of its components.  A similar relation to homology is contained in
Proposition~\ref{truncations:non-negative}(2).

  \begin{sublemma}
\label{levels:estimates}
For every complex of $R$-modules $M$ there is an inequality
  \[
\level{R}{\cc}M \leq
\inf\bigg\{\sum_{n\in\BZ}\level{R}{\cc}{L_n}\ \left|\
L\simeq M\text{ in }\dcat R\bigg\}\right.
  \]
   \end{sublemma}

  \begin{proof}
As $\level R{\cc}M=\level R{\cc}L$, it suffices to assume that only
finitely many components $L_{n}$ are non-zero and to show that then one has
  \[
\level{R}{\cc}L \leq\sum_{n\in\BZ}\level{R}{\cc}{L_n}\,.
  \]
The complex $L$ admits a filtration $\cdots \subseteq L_{\les n-1}
\subseteq L_{\les n}\subseteq\cdots$, and there are isomorphisms $L_{\les
n}/L_{\les n-1}\cong \shift^nL_n$. Lemma~\ref {levels:properties}(1)
yields $\level{R}{\cc}{\shift^nL_n}= \level R{\cc}{L_n}$, so one
gets the desired inequality by a repeated application of Lemma
\ref{levels:properties}(2).
  \end{proof}

As an example, we deduce an inequality that appears in the introduction.

  \begin{subexample}
  \label{plevel:card}
If $P$ is a complex of finite projective $R$-modules, then
  \[
\level{R}{R}P\leq\card\{n\in\BZ\mid P_n\ne0\}\,.
  \]

Just note that each $P_n$ is in $\thickn1RR$ and apply the
Lemma~\ref{levels:estimates}.
   \end{subexample}
    \end{bfchunk}

\section{Levels of DG modules}
\label{Levels of DG modules}

In this section we move from modules over rings to DG (= differential
graded) modules over DG algebras.  We start by recalling some concepts
and collecting notation and basic facts concerning DG algebras and
DG modules.

Throughout this paper, we consider DG objects as collections of abelian
groups indexed by the integers, rather than direct sums of such groups.
This point view, prevalent among topologists, is systematically developed
in MacLane's books; see especially \cite[\S{VI.2}]{Ml}.  A consequence
is that every element $m$ of a graded object $M$ has a degree, denoted
$|m|$; namely, one has $|m|=i$ if and only if $m$ belongs to $M_i$.
Differentials have degree $-1$.  DG algebras act on their DG
modules from the left.

Rings are identified with DG algebras concentrated in degree $0$.
DG modules over a ring are just complexes of modules, and modules over
it are identified with complexes concentrated in degree $0$.

A \emph{graded algebra} is a DG algebra with zero differential.  It should
be noted that in a DG algebra elements of different degrees cannot be
added.  Furthermore, DG modules over a graded algebra $A$ should not be
confused with complexes over it:  The first ones have a unique grading
and their differentials satisfy the identity $\dd(am)=(-1)^{|a|}a\dd(m)$,
while the second ones are equipped with two gradings (homological and
internal), and their differentials are $A$-linear.

Let $A$ be a DG algebra and $M$, $N$ be left DG modules.  

A \emph{homomorphism} $\beta\col M\to N$ of degree $d$
is a family $(\beta_i\col M_i\to N_{i+d})_{i\in\BZ}$ of
additive maps satisfying $\beta(am)=(-1)^{d|a|}a\beta(m)$ for
$a\in A$ and $m\in M$.  All such homomorphisms form the $d$th
component of a complex, $\Hom AMN$, whose differential is given by
$\dd(\beta)=\dd^N\beta-(-1)^{|\beta|}\beta\dd^M$. A \emph{morphism} of
DG modules $M\to N$ is a homomorphism $\beta$ of degree $0$, satisfying
$\dd^N\beta=\beta\dd^M$.  A \emph{quasi-isomorphism} is a morphism that
induces an isomorphism in homology.

We let $A^\natural$ and $M^\natural$ denote the underlying
graded algebra and graded module over it.  For each integer
$s$ let $\shift^sM^\natural$ denote the graded $\BZ$-module
with $(\shift^sM^{\natural})_i=M_{i-s}$.  The map $m\mapsto m$ is a
homomorphism $\sigma^s\col M^{\natural}\to\shift^sM^\natural$ of degree
$s$.  Let $\shift^{s} M$ denote the DG $A$-module with underlying
graded $\BZ$-module $\shift^sM^\natural$, action of $A$ given
by $a\sigma^s(m)=(-1)^{|a|s}\sigma^s(am)$, and differential by
$\dd(\sigma^s(m))=(-1)^s\sigma^s(\dd(m))$.

Proofs of statements in \ref{semi-free}--\ref{resolutions:functors}
below can be found in \cite{AFH}.

\begin{bfchunk}{Semi-free DG modules.}
\label{semi-free} 
A DG module $F$ over a DG algebra $A$ is \emph{semi-free} if it admits a family
$\Filt Fn_{n\in\BZ}$ of DG $A$-submodules satisfying the conditions:
\begin{gather*}
\filt Fn \subseteq \filt F{n+1}\,, \quad \filt F{-1}=0\,,
\quad \bigcup_{n\in\BZ}\filt Fn=F\,,\quad\text{and}
\\
\filt F{n+1}/\filt F{n} \quad
\text{is isomorphic to a direct sum of suspensions of }A\,.
  \end{gather*}
In this case, the functors $\Hom AF-$ and $-\otimes_AF$, defined on
the category of DG $A$-modules and right DG $A$-modules, respectively,
preserve quasi-isomorphisms.

Thus, if $A$ is an algebra over a field $k$ the functor $\Hom A{-}{\Hom
kFk}$ preserves quasi-isomorphisms of right DG $A$-modules, as it is
isomorphic to $\Hom k{-\otimes_{A}F}k$.
  \end{bfchunk}

  \begin{bfchunk}{Semi-free resolutions.}
 \label{resolutions:semi-free}
Each DG $A$-module $M$ admits a quasi-isomorphism $F\to M$ with $F$
a semi-free DG $A$-module.  Such a \emph{semi-free resolution} of $M$
is unique up to homotopy of DG $A$-modules.
 \end{bfchunk}

 \begin{bfchunk}{Derived categories.}
Let $A$ be a DG algebra.  We let $\dcat A$ denote the derived category of DG $A$-modules. Its objects are DG
$A$-modules, and it can be realized as the homotopy category of semi-free DG $A$-modules.
The derived category is triangulated, see \cite[\S4]{Ke} or \cite[3.2]{Kr}
for constructions.

We consider every ring $R$ as a DG algebra concentrated in degree zero.
Its DG modules are then simply the complexes of $R$-modules, and the
derived category $\dcat R$ coincides with the derived category of $R$-modules.
\end{bfchunk}

\begin{bfchunk}{Derived functors.}
 \label{resolutions:functors}
For each DG $A$-module one sets
 \[
\Rhom AM- = \Hom AF-
\quad\text{and}\quad
(-\lotimes AM)=(-\otimes_AF)\,,
 \]
where $F\to M$ is some semi-free resolution.  This yields well defined
exact functors on the derived category of DG $A$-modules and right DG
$A$-modules, respectively.
  \end{bfchunk}

Given a subcategory $\cc$ of $\dcat A$ and a DG $A$-module $M$ we write
$\thickn nA{\cc}$ and $\level A{\cc}M$ in place of $\thickn n{\dcat
A}{\cc}$ and $\level {\dcat A}{\cc}M$, respectively.
 
We record some easy consequences of general properties of levels.

  \begin{chunk}
Let $A$ be a DG algebra, $M$ a DG $A$-module, and
$\cc\subseteq\dcat A$ a subcategory. 

 \begin{subchunk}
\label{levels:finiteness1}
If $M$ has a filtration by DG submodules
    \[
0=\filt{M}{0} \subseteq \filt{M}{1}\subseteq \cdots \subseteq \filt{M}{m}
= M
    \]
then the following inequality holds:
    \[
\level{A}{\cc}M \leq
\sum_{i=1}^m\level{A}{\cc}{\filt{M}{i}/\filt{M}{i-1}}\,.
    \]

Indeed, this is seen through iterated applications of Lemma~\ref
{levels:properties}(2) to the triangles in $\dcat A$ defined by the
exact sequences of DG $A$-modules
  \begin{equation*}
0\to\filt{M}{i-1}\to\filt{M}{i} \to \filt{M}{i}/\filt{M}{i-1}\to0 
  \end{equation*}
 \end{subchunk}

 \begin{subchunk}
\label{levels:finiteness2}
Let $B$ be a DG algebra and $L$ a left-right $B$-$A$-bimodule with
$\hh{L\lotimes AC}$ a noetherian, respectively, artinian, $\hh B$-module for each 
$C\in\cc$.  If $\level A{\cc}M$ is finite, then the $\hh B$-module 
$\hh{L\lotimes AM}$ is noetherian, respectively, artinian.
  \end{subchunk} 

Indeed, as $(L\lotimes A-)$ is an exact functor from $\dcat A$ to $\dcat B$,
Lemma~\ref{levels:properties}(6) shows that $L\lotimes AM$ has finite
$(L\lotimes A\cc)$-level, so the desired assertion follows from 
\ref{levels:inheritance}.
    \end{chunk}

Next we consider how levels behave under changes of DG algebras.

  \begin{bfchunk}{Morphisms of DG algebras.}
Let $\vf\col A\to B$ be a morphism of DG algebras.

The following property is used implicitly in many arguments.

  \begin{subchunk}
     \label{morphisms:adjoint}
The map $\vf$  induces an adjoint pair of functors of triangulated categories
    \[
\xymatrixcolsep{4.5pc}
\xymatrix{
\dcat A\ar@{->}[r]<.6ex>^-{(B\lotimes A-)}
&\dcat B\ar@{->}[l]<.6ex>^-{\vf_*}
}
    \]
where $\vf_*$ is the functor restricting the action of $B$ to $A$.
In particular, for every $M\in\dcat A$ and every $N\in\dcat B$ there
are canonical morphisms
    \[
M\lra \vf_*(B\lotimes AM) \qquad\text{and}\qquad B\lotimes A\vf_*(N)\lra N
    \]
  \end{subchunk}

   \begin{subchunk}
     \label{morphisms:equivalence}
If $\vf$ is a quasi-isomorphism, then $(B\lotimes A-)$ and $\vf_*$
are inverse equivalences.
    \end{subchunk}

Indeed, with $U$ a semi-free resolution of the DG module $M$ over $A$, see
\ref{resolutions:semi-free}, the morphism $M\to \vf_*(B\lotimes AM)$ is
represented by $\vf \otimes_{A}U$; it is a quasi-isomorphism since $\vf$
is one.  Let $\nu\col V\to N$ be a semi-free resolution of $\vf_{*}(N)$
over $A$. The morphism $B\lotimes A\vf_*(N)\to N$ is represented by
the morphism $\mu\col B\otimes_{A}V\to N$, where $b\otimes v\mapsto
b\nu(v)$. The quasi-isomorphism $\nu$ factors as
 \[
V\xra{\vf\otimes_{A}V} B\otimes_{A}V \xra{\mu} N\,.
 \]
Since $V$ is semi-free, $\vf\otimes_{A}V$ is a quasi-isomorphism, and
hence so is $\mu$.

  \begin{subchunk}
\label{morphisms:objects}
Let $M$ and $N$ be DG modules over $A$ and $B$, respectively, and let
$\mu\col M\to N$ be a morphism of complexes of abelian groups satisfying
$\mu(am)=\vf(a)\mu(m)$ for all $a\in A$ and $m\in M$.  If $\vf$ and $ \mu$
are quasi-isomorphisms, then in $\dcat A$ and $\dcat B$, respectively,
one has canonical isomorphisms
   \[
M\simeq\vf_*(N)
    \qquad\text{and}\qquad
B\lotimes AM \simeq N
   \]
Indeed, $M\simeq\vf_*(N)$ holds by assumption, so
\ref{morphisms:equivalence} yields isomorphisms
\[
B\lotimes AM\simeq B\lotimes A\vf_*(N)\simeq N\,.
\]
\end{subchunk}
      \end{bfchunk}

Now we can track the behavior of levels under change of DG algebras.

  \begin{proposition}
\label{levels:basechange2}
Let $\vf\col A\to B$ be a morphism of DG algebras.

For all DG $A$-modules $C,M$ and all DG $B$-modules $D,N$ the
following hold.
  \begin{enumerate}[\quad\rm(1)]
    \item
There are inequalities
  \begin{align*}
\level{A}{C}M &\ge\level{B}{B\lotimes AC}{B\lotimes AM}
\\
\level{A}{\vf_*(D)}{\vf_*(N)}&\le\level{B}{D}N
  \end{align*}
Equalities hold when $\vf$ is a quasi-isomorphism.
    \item
If both $\level{A}{A}{\vf_*(B)}$ and $\level{B}BN$ are finite, then
so is  $\level{A}{A}{\vf_*(N)}$.
  \end{enumerate}
    \end{proposition}

  \begin{proof}
(1) Lemma~\ref{levels:properties}(6) yields the inequalities and shows
that they become equalities if $B\lotimes A-$ and $\vf_*$ are
equivalences; now refer to \ref{morphisms:equivalence}.

(2) Since $\level{B}BN$ is finite, (1) shows that so is
$\level{A}{\vf_*( B)}{\vf_*(N)}$.  It follows that
$\vf_*(N)$ is in $\thick{A} {\vf_*(B)}$.
As $\level{A}A{\vf_*(B)}$ is finite, one has an inclusion
  \[
\thick{A}{\vf_*(B)}\subseteq \thick{A}A
  \]
Thus, $\vf_*(N)$ is in $\thick{A}A$, which means  $\level{A}A{\vf_*(N)}$
is finite, as desired.
  \end{proof}

The balance of this section deals with special properties of DG algebras
concentrated either in non-negative degrees or in non-positive degrees.

  \begin{chunk}
  \label{augmentation:non-negative}
A DG algebra $A$ is \emph{non-negative} if it has $A_n=0$ for all $n<0$.

When $A$ is a non-negative DG algebra the subcomplex
\[
J=\quad \cdots\lra A_2\lra A_1\lra\Image(\dd_1)\lra0\lra\cdots
\]
is a DG ideal of $A$, and $A/J$ is naturally isomorphic to $\HH 0A$.  The
morphism of DG algebras $\eps\col A\to\HH0A$ is called the \emph{canonical
augmentation} of $A$.  Via the functor $\eps_*\col\dcat{\HH0A}\to\dcat A$
we identify complexes of $\HH 0A$-modules with DG $A$-modules; the same
letter denotes a complex in $\dcat{\HH0A}$ and its image in $\dcat A$.
  \end{chunk}

  \begin{proposition}
\label{dg:h0actions}
Let $A$ and $B$ be non-negative DG algebras.

If $\SF\col\dcat A\to\dcat B$ is an equivalence of triangulated categories
induced by a chain of quasi-isomorphisms of DG algebras, then in $\dcat B$
one has an isomorphism
    \[
\SF(\HH0A)\simeq\HH0B
    \]
  \end{proposition}

  \begin{proof}
By hypothesis, there exists a sequence of quasi-isomorphisms
    \[
  \xymatrixcolsep{2.35pc}
  \xymatrix{
A\ar@{->}[r]^-{\simeq}
&A^0\ar@{<-}[r]^-{\simeq}
&A^1\ar@{->}[r]^-{\simeq}
&\ \cdots\ \ar@{<-}[r]^-{\simeq}
&A^{i-1}\ar@{->}[r]^-{\simeq}
&A^{i}\ar@{<-}[r]^-{\simeq}
&B\
}
  \]
of DG algebras.  For $j=0,\dots,i$ the subcomplex
  \[
A_+^j=\quad \cdots\lra A^j_2\lra A^j_{1}\lra
\Ker\big(\dd_0^{A^j}\big)\lra0\lra\cdots
  \]
is a functorially defined DG subalgebra of $A^j$, and the inclusion
$A_+^j\subseteq A^j$ is a quasi-isomorphism.  Thus, the original
sequence induces a commutative diagram
    \[
  \xymatrixcolsep{1pc}
  \xymatrixrowsep{2pc}
\xymatrix{
A\ar@{->}[r]^-{\simeq}\ar@{->}[d]
&A_+^0\ar@{<-}[r]^-{\simeq}\ar@{->}[d]
&A_+^1\ar@{->}[r]^-{\simeq}\ar@{->}[d]
&\ \cdots\ \ar@{<-}[r]^-{\simeq}
&A_+^{i-1}\ar@{->}[r]^-{\simeq}\ar@{->}[d]
&A_+^{i}\ar@{<-}[r]^-{\simeq}\ar@{->}[d]
&B\ar@{->}[d]
  \\
\HH0{A}\ar@{->}[r]^-{\cong}
&\HH0{A_+^0}\ar@{<-}[r]^-{\cong}
&\HH0{A_+^1}\ar@{->}[r]^-{\cong}
&\ \cdots\ \ar@{<-}[r]^-{\cong}
&\HH0{A_+^{i-1}}\ar@{->}[r]^-{\cong}
&\HH0{A_+^{i}}\ar@{<-}[r]^-{\cong}
&\HH0{B}
}
  \]
of isomorphisms of DG algebras, where all DG algebras in the top row
are non-negative and all vertical arrows are canonical augmentations.
Iterated applications of \ref{morphisms:objects} produce the desired
isomorphism in $\dcat B$.
   \end{proof}

  \begin{proposition}
\label{truncations:non-negative}
Let $A$ be a non-negative DG algebra.

For every DG $A$-module the following hold.
    \begin{enumerate}[\quad\rm(1)]
    \item
If $\inf\{n\in\BZ\var \HH nM\ne0\}=i>-\infty$, then
there is an exact triangle
\[
M' \lra M\lra \shift^i(\HH iM)\lra \shift M'
\]
in $\dcat A$, with $\HH i{M'}=0$ and $\HH n{M'}\cong \HH n{M}$
for $n\ne i$.
    \item
For every class $\cc$ of $\HH 0A$-modules one has
\[
\level{A}{\cc}M \leq \sum_{n\in\BZ} \level{\HH 0A}{\cc}
{\HH nM}\,.
\]
  \end{enumerate}
    \end{proposition}

  \begin{proof}
(1)  Since $A$ is non-negative, the following subcomplexes
    \[
\xymatrixcolsep{1.5pc}
\xymatrixrowsep{.7pc}
\xymatrix{
M' = \quad\cdots\ar@{->}[r]
&M_{i+1}\ar@{->}[r]
&\Image(\dd_{i+1})\ar@{->}[r]
&0\ar@{->}[r]
&\cdots
\\
M''= \quad\cdots\ar@{->}[r]
&M_{i+1}\ar@{->}[r]
&\Ker(\dd_i)\ar@{->}[r]
&0\ar@{->}[r]
&\cdots
}
    \]
of $M$ are closed under multiplication by elements on $A$; in other
words, they are DG submodules of $M$. It remains to observe that the
inclusion $M''\subseteq M$ is a quasi-isomorphism, and that one has an
exact sequence of DG $A$-modules
    \[
0\lra M' \lra M'' \lra \shift^i\HH iM\lra 0
    \]
It gives rise to an exact triangle with the desired properties.

(2) We may assume that the number $w(M)=\card\{n\in\BN\var \HH nM\ne0\}$
is finite.  As $w(M)=0$ means $M\simeq 0$, the desired inequality is
evident in this case, so we may assume that it holds for DG modules $N$
with $w(N)<r$ for some $r\ge1$.

Set $i=\inf\{n\var \HH nM\ne0\}$.  Lemma~\ref{levels:properties}(2)
applied to the exact triangle in (1) and Lemma~\ref
{levels:properties}(1)
yield the first relation below:
    \begin{align*}
\level{A}{\cc}M
&\leq\level{A}{\cc}{\HH iM}
    +\level{A}{\cc}{M'}
\\
&\leq\level{\HH 0A}{\cc}{\HH iM}+\level{A}{\cc}{M'}
\\
&\leq\level{\HH 0A}{\cc}{\HH iM}+\sum_{n\ges i+1}\level{\HH 0A}{\cc}
{\HH n{M'}}
\\
&=\sum_{n\ges i}\level{\HH 0A}{\cc}{\HH n{M}}\,.
  \end{align*}
Proposition~\ref{levels:basechange2}(1), applied to the morphism $\eps\col A\to \HH 0A$,
yields the second one.  As $w(M')=r-1$ by part (1), the third inequality
is the induction hypothesis.  The equality uses the expressions for the
modules ${\HH n{M'}}$, again from (1).
   \end{proof}

  \begin{chunk} \label{augmentation:non-positive}
A DG algebra $A$ with $A_n=0$ for $n>0$ is said to be \emph{non-positive}.
Such a DG algebra has no augmentation to $\HH0A$ in general.  However,
the subcomplex
   \[
J=\quad \cdots\lra 0\lra A_{-1}\lra A_{-2}\lra\cdots
   \]
is a DG ideal of $A$. Thus, when $\dd(A_0)=0$ one has natural isomorphisms
$A/J\cong A_0\cong\HH 0A$, and hence a \emph{canonical augmentation}
$\eps\col A\to \HH 0A$.  As in \ref{augmentation:non-negative},
for every complex $M\in\dcat{\HH0A}$ we let $M$ denote also the DG
$A$-module $\eps_*(M)$.
  \end{chunk}

  \begin{proposition}
\label{truncations:non-positive}
If $A$ is a non-positive DG algebra with $\dd(A_0)=0$ and
$A_0$ semi-simple, and $M$ is a DG $A$-module, then the following hold.
    \begin{enumerate}[\quad\rm(1)]
    \item
If $\sup\{n\in\BZ\mid \HH nM\ne0\}=s<\infty$, then there is an exact
triangle
  \[
M' \lra M\lra \shift^s\HH sM\lra \shift M'
  \]
in $\dcat A$, with $\HH s{M'}=0$ and $\HH n{M'}\cong\HH n{M}$
for $n\ne s$.
  \item
For every class $\cc$ of $\HH 0A$-modules  one has
    \[
\level{A}{\cc}M \leq \sum_{n\in\BZ} \level{\HH 0A}{\cc}
{\HH nM}\,.
  \]
   \end{enumerate}
    \end{proposition}

  \begin{proof}
(1) The hypothesis $\dd(A_0)=0$ implies $\dd^M$ is $A_0$-linear. As
$A_0$ is semi-simple,
    \begin{gather*}
0\lra \Ker(\dd_s)\lra M_s\lra \Image(\dd_s)\lra 0\\ 0\lra
\Image(\dd_{s+1}) \lra \Ker(\dd_s)\lra \HH sM \lra 0
    \end{gather*}
are split-exact sequences of $A_0$-modules.  Thus, there are
isomorphisms
    \[
M_s\cong\Image(\dd_{s+1})\oplus\HH sM\oplus C
\quad\text{and}\quad
C\cong \Image(\dd_s)
    \]
of $A_0$-modules, the latter induced by $\dd_s$.  As $A$ is
non-positive,
    \[
\xymatrixcolsep{1.5pc}
\xymatrixrowsep{.7pc}
\xymatrix{
M' = \quad\cdots\ar@{->}[r]
&0\ar@{->}[r]
&C\ar@{->}[r]^-{\dd_s}
&M_{s-1}\ar@{->}[r]
&\cdots
\\
M''= \quad\cdots\ar@{->}[r]
&0\ar@{->}[r]
&\HH sM\oplus C\ar@{->}[r]^-{\dd_s}
&M_{s-1}\ar@{->}[r]
&\cdots
}
    \]
are DG $A$-submodules of $M$ and there is an exact sequence of DG
$A$-modules
    \[
0\lra M' \lra M''\lra \shift^s\HH sM\lra 0
    \]
The inclusion $M''\subseteq M$ is a quasi-isomorphism, so it yields
the desired exact triangle; the homology of $M'$ is computed from the
homology exact sequence.

(2) follows from (1) by an argument parallel to that for Proposition
\ref{truncations:non-negative}(2).
     \end{proof}

\section{Perfect DG modules}
\label{Perfect DG modules}

We say that a DG module $M$ over a DG algebra $A$ is \emph{perfect}
if $\level AAM$ is finite.  The first result of this section describes
the structure of perfect DG modules.  Specialized to the case of rings
it shows that our terminology is consistent with the traditional notion
for complexes, see \ref{Levels of complexes of modules}.  We extend this
characterization to DG modules over non-negative DG algebras.

Over certain DG algebras which arise in many applications we establish
a homological and hence verifiable criterion for perfection. We finish
the section with examples showing that the hypotheses of this last result
cannot be relaxed easily.

  \begin{bfchunk}{Semi-freeness.}
    \label{plevel:semifreeness}
A \emph{semi-free filtration} of a DG $A$-module $F$ is a family
$\Filt Fn_{n\in\BZ}$ of DG $A$-submodules satisfying the conditions in
\ref{semi-free}. Such a filtration $\Filt Fn_{n\in\BZ}$ has \emph{class
at most} $l$ if $F^l=F$ holds for some integer $l$; it is \emph{finite}
if, in addition, its subquotients are finitely generated.

A DG $A$-module admitting a (finite) semi-free filtration (of class at
most $l$) is said to be (\emph{finite}) \emph{semi-free} (\emph{of class
at most $l$}).  Note that $0$ is the only DG module that is semi-free
of class at most $-1$.
  \end{bfchunk}

The next theorem suggests that the $A$-level of a DG $A$-module may be
thought of as a kind of `projective dimension'.

  \begin{theorem}
\label{plevel:filtrations}
Let $A$ be a DG algebra and $l$ a non-negative integer.

A DG $A$-module $M$ has $\level AAM\leq l$ if and only if it is a retract
of some finite semi-free DG module of class at most $l-1$.
  \end{theorem}

  \begin{proof}
Let $\cf$ denote the full subcategory $\sadd{A}$ of $\dcat A$; its objects
are the DG modules isomorphic to finite direct sums of suspended copies of
$A$, see \ref{Additive closures}.  In view of \ref{layers:symmetry},
it suffices to prove that $\cf^{\star l}$ is the smallest strict
subcategory of $\dcat A$ that contains all finite semi-free DG modules
of class at most $l-1$.

Indeed, a semi-free filtration $\Filt Fn_{n\in\BZ}$ of class at most
$l-1$
yields exact triangles
    \[
\filt{F}{n-1} \lra \filt{F}{n} \lra \filt{F}{n}/\filt{F}{n-1}\lra \shift
\filt{F}{n-1}
     \quad\text{for}\quad 1\leq  n\leq l-1\,.
    \]
Since $\filt{F}{0}$ and $\filt{F}{n}/\filt{F}{n-1}$ are in $\cf$,
induction shows that $F$ is in $\cf^{\star l}$.

Conversely, assume that $M$ is in $\cf^{\star l}$ for some $l\geq0$. By induction on $l$ we prove that $M$ is isomorphic to a finite semi-free DG module $F$ of class at most $l-1$.  For $l=0$ the assertion is evident. For $l\geq 1 $ one has in $\dcat A$ an exact triangle
    \[
G\xra{\ \gamma\ }L\lra M\lra \shift G
    \]
with $L$ in $\cf^{(l-1)\star}$ and $G$ in $\cf$.  By the induction
hypothesis, $L$ is isomorphic to a DG $A$-module with a finite semi-free
filtration of class $l-2$.  Thus, there is an exact triangle as above,
where $L$ has a finite semi-free filtration $\Filt Ln$ with $L^{l-2}=L$.

Since $G$ is in $\cf$, the triangle above is isomorphic to a triangle
    \[
\wt G\xra{\ \wt\gamma\ }L\lra M\lra \shift\wt G
    \]
where $\wt G$ is a finite direct sum of suspended copies of $A$.
Set $F=\cone{\wt\gamma}$.  One then has $M\simeq F$ and an exact sequence
of DG $A$-modules
    \[
0  \lra L \xra{\ \lambda\ } F \lra \shift  \wt G\lra 0\,.
    \]
One gets a finite semi-free filtration of $F$ with $\filt F{l-1}=F$ 
by setting
    \[
\filt Fn =
    \begin{cases}
\lambda(\filt Ln) & \text{for }n\leq l-2\,; \\
         F  & \text{for } n\geq l-1\,.
    \end{cases} \]
It shows that the class of $F$ is at most $l-1$, as desired.
     \end{proof}

By Theorem~\ref{plevel:filtrations} a DG $A$-module $M$ of finite
$A$-level is isomorphic in $\dcat A$ to a DG module $P$ with $P^\natural$
finite projective over $A^\natural$.  Next we prove that the converse
holds when $A$ is non-negative, in particular, when $M$ is a complex
over a ring.

In the argument, we use the following classical fact on the structure
of graded projective modules over graded rings; see \cite[6.6]{Sw}.

\begin{chunk}
\label{level:extended}
Let $B$ be a non-negatively (respectively, non-positively) graded 
algebra.

For each bounded below (respectively, bounded above) projective graded
$B$-module $N$ there is an isomorphism of graded $B$-modules $N\cong
V\otimes_{B_0}B$, with $V$ a bounded below (respectively, bounded above)
projective graded $B_0$-module.

The $B_0$-module $V$ is defined uniquely up to isomorphism: $V\cong
(B/J)\otimes_BN$, where $J$ denotes the ideal of elements of positive
(respectively, negative) degree.
  \end{chunk}

  \begin{proposition}
\label{level:semi-projective}
Let $A$ be a non-negative DG algebra.

A bounded below DG $A$-module $M$ is a direct summand of some
(finite) semi-free DG $A$-module if and only if the underlying graded
$A^\natural$-module is (finite) projective.
  \end{proposition}

  \begin{proof}
The `only if' part is evident.  For the converse, set $R=A_0$.  By
\ref{level:extended}, one has $M^\natural\cong V\otimes_RA^{\natural}$
for some bounded below graded $R$-module $V$ with each $V_n$ projective.
Pick for each $n\in\BZ$ an $R$-module $W_n$ so that $V_n\oplus W_n$
is free.  If $M^\natural$ is finite over $A^\natural$, each $V_n$ is
finite over $R$, so choose $W_n$ finite, as one may.  Let $W$ be the
complex of $R$-modules having $W_n$ as $n$th component and $\dd^W=0$.

Form the DG $A$-module $F=M\oplus(W\otimes_RA)$.  As $A$ is non-negative,
  \[
F^n=(V\oplus W)_{\les n}\otimes_RA
  \]
is a DG submodule of $F$, and
$\cdots \subseteq F^n \subseteq F^{n+1}\subseteq\cdots$ is a semi-free
filtration of $F$; it is finite when the $R$-module $W$ is finite.
  \end{proof}

Theorem \ref{plevel:filtrations} and Proposition \ref{level:semi-projective}
yield:

  \begin{corollary}
\label{level:DGperfection1}
A DG module over a non-negative DG algebra $A$ is perfect if it is
quasi-isomorphic to a DG module $P$ with $P^\natural$ finite projective
over $A^\natural$.
 \qed
  \end{corollary}

{F}rom the preceding corollary and Proposition~\ref{levels:basechange2}(2) 
one gets:

  \begin{corollary}
\label{levels:basechange3}
Let $\vf\col A\to B$ be a morphism of DG algebras such that $A$ is
non-negative and the graded $A^\natural$-module $\vf_*(B)^\natural$
is finite projective.

If $N$ is a perfect DG $B$-module, then the DG $A$-module $\vf_*(N)$ 
is perfect.
 \qed
    \end{corollary}

Our next goal is to obtain a homological criterion for perfection.
We start by noting that \ref{levels:finiteness2} applied to the class $\cc=\{A\}$ gives
a homological obstruction:

 \begin{remark}
  \label{plevel:lotimes}
Let $A$ and $B$ be DG algebras, $L$ a left-right $B$-$A$-bimodule,
and $M$ a perfect  DG $A$-module.

If the graded $\hh B$-module $\hh L$ is noetherian (respectively, artinian) then
the graded $\hh B$-module $\hh{L\lotimes AM}$ is noetherian (respectively, artinian).
    \end{remark}

Under additional hypotheses, we prove that finiteness is
the only obstruction; for the definition of the canonical
augmentation, see Remarks~\ref{augmentation:non-negative} and
\ref{augmentation:non-positive}.

  \begin{theorem}
  \label{level:DGperfection2}
Let $A$ be a DG algebra and $M$ a DG $A$-module, such that
  \begin{enumerate}[\quad\rm(a)]
 \item
$A$ is non-negative, $\dd(A_1)=0$, and $\hh M$ is bounded below; or
 \item
$A$ is non-positive, $A_{-1}=0$, and $\hh M$ is bounded above.
  \end{enumerate}
Set $k=A_0$ and let $\eps\col A\to\HH 0A=k$ denote the canonical 
augmentation.

If $k$ is a field, then the following conditions are equivalent.
\begin{enumerate}[\quad\rm(i)]
 \item
$M$ is perfect over $A$.
 \item
$\hh{k\lotimes AM}$ is finite over $k$.
 \item
$M\simeq F$ for some finite semi-free $F\in\dcat A$.
 \item
$M\simeq F$ for some $F\in\dcat A$ with $F^\natural$ a finite projective
graded $A^\natural$-module.
  \end{enumerate}
When they are satisfied, the inequality below holds:
 \[
\level AAM\le\card\{n\in\BZ\mid\HH n{k\lotimes AM}\ne0\}\,.
 \]
   \end{theorem}

In our argument, we use the existence of minimal semi-free resolutions,
see \cite{AFH}:

\begin{bfchunk}{Minimal resolutions.}
 \label{resolutions:minimal}
Let $A$ be a DG algebra, and $M$ a DG $A$-module, as in
Theorem~\ref{level:DGperfection2}.  Assume $\HH 0A$ is a field. Set
$J=\{a\in A : |a|\ne0\}$; this is a DG ideal of $A$, and there exists
a semi-free resolution $F\to M$ with the property $\dd(F)\subseteq JF$.
 \end{bfchunk}

  \begin{proof}[Proof of Theorem~\emph{\ref{level:DGperfection2}}]
Set $J=\Ker(\eps)$.

(i) $\implies$ (iv).  This follows from Theorem~\ref{plevel:filtrations}.

(iv) $\implies$ (iii).  By \ref{level:extended} the $A^\natural$-module
$F^\natural$ is free.  We induce on $r=\rank_{A^\natural}F^{\natural}$.

For $r=0$ the assertion is obvious.  Fix $r\ge1$ and assume that
the assertion holds for all DG modules with underlying graded module
of smaller rank.  Choose a non-zero element $e$ of $F$, with $|e|$
minimal in case (a) and maximal in case (b).  One can then find a basis
$\{e_1,\dots,e_r\}$ of $F^\natural$ over $A^\natural$, with $e_1=e$.

If $\dd(e)=0$, then $\ov F=F/Ae$ is a DG $A$-module with $(\ov
F)^{\natural}$ free of rank $r-1$, so by the induction hypothesis it
has a finite  semi-free filtration $\Filt{\ov F}n$.  Setting $F^0=Ae$
and letting $F^n$ denote the inverse image in $F$ of $\ov F{}^{n-1}$ for
$n\ge 1$, we obtain a finite semi-free filtration of $F$.  This covers
case (a), as well as case (b) for $r=1$.

Finally, assume (b) holds, $r\ge2$, and $\dd(e)=f\ne0$.  Set $|e|=j$
and $J=A_{\le-2}$.  As $|f|=j-1$, and $(JF)_i=0$ for $i\ge j-1$,
we get $f\notin JF$.  By Nakayama's Lemma $F^\natural$ has a  basis
$\{e_1,\dots,e_r\}$ with $e_1=e$ and $e_2=f$.  The graded submodule
$E^\natural$ of $F$ generated by $e$ and $f$ is a DG submodule, and has
$\hh E=0$.  The exact sequence of DG modules $0\to E\to F\to F/E\to 0$
yields a quasi-isomorphism $F\simeq F/E$ and shows that $(F/E)^\natural$
is free of rank $r-2$ over $A^\natural$.  Thus we get $M\simeq F/E$ in
$\dcat A$, and the induction hypothesis implies that $F/E$ is semi-free.

(iii) $\implies$ (ii).  This is due to the isomorphism $\hh{k\lotimes
AM}\cong\hh{k\otimes_AF}$.

(ii) $\implies$ (i).  
It suffices to show that (ii) implies the inequality in the statement
of the theorem.  Setting $V=\hh{k\lotimes AM}$ we argue by induction
on $v=\rank_kV$.  

Choose a quasi-isomorphism $F\to M$ with $F$ a semi-free DG $A$-module
and $\dd(F)\subseteq JF$; see \ref{resolutions:minimal}.  One then
has isomorphisms of graded $k$-vector spaces
 \[
k\otimes_AF=\hh{k\otimes_AF}\cong\hh{k\lotimes AM}=V
 \]
{F}rom \ref{level:extended} one gets $F^\natural\cong
V\otimes_kA^\natural$.  In particular, $v=0$ implies $F=0$, and hence
$\level AAF=0$.  This is the basis for our induction.

Assume $v\ge1$.  For $i=\inf\{n\in\BZ\var V_n\ne0\}$ and
$s=\sup\{n\in\BZ\var V_n\ne0\}$ we have
  \begin{alignat*}{2}
\dd(V_i\otimes_k 1)\subseteq (V\otimes_kA_{\ges1})_{i-1}=0 &&\quad\text{in
case (a);}
 \\
\dd(V_s\otimes_k 1)\subseteq
(V\otimes_kA_{\les-1})_{s-1}=V_s\otimes_kA_{-1}=0 &&\text{in case (b).}
 \end{alignat*}
Setting $j=i$ in case (a) and $j=s$ in case (b), we get an exact sequence
 \begin{gather*}
0\lra(\shift^jV_j)\otimes_kA\lra F\lra G\lra 0
 \end{gather*}
of DG $A$-modules.  The homology exact sequence of the functor $(k\lotimes
A-)$ splits:
 \[
0\lra\shift^jV_j\lra V\lra \hh{k\lotimes AG}\lra 0
 \]
In particular, we get $\rank_k\hh{k\lotimes AG}<v$.  The DG module $G$
is bounded below in case (a) and above in case (b), so the induction
hypothesis applies to it.  Since $\level AA{V\otimes_kA}=1$ holds by
definition, \ref{levels:finiteness1} yields the inequality below:
 \begin{align*}
\level AAM
&=\level AAF\\
&\le 1+\level AAG\\
&=1+\card\{n\in\BZ\mid\HH n{k\lotimes AG}\ne0\}\\
&=\card\{n\in\BZ\mid\HH n{k\lotimes AM}\ne0\}
 \end{align*}
The last equality comes from the homology exact sequence above.
  \end{proof}

  \begin{remark}
 \label{level:others}
The proof above depends on the existence of minimal semi-free resolutions, which are more widely available.  For example, when $A$ is non-negative and $M$ is bounded below they exist if the ring $A_0$ is artinian, and
also if $A_0$ is local, $\hh A$ is noetherian, and $\hh M$ is finite over $\hh A$.  Minor modifications in our arguments extend to these cases the validity of the theorem. 

In the sequel we use only the following special case.
 \end{remark}

\begin{proposition}
\label{level:DGperfection3}
Let $A$ be a DG algebra with zero differential and let $M$ be a DG
$A$-module, such that  there is an isomorphism $M\simeq\hh M$ in
$\dcat A$.

If $A$ is non-negative or non-positive, and the ring $A_{0}$
is Artinian and local, then $M$ is perfect if and only if the
graded $A^\natural$-module $\hh M$ has a finite free resolution.
\end{proposition}

 \begin{proof}
Let $k$ denote the residue field of $A_{0}$ and $J$ the kernel of the
canonical augmentation $A\to k$. The graded $A^\natural$-module $N=\hh M$
has a \emph{minimal} free resolution
   \[
F\vphantom{I}^{\scriptscriptstyle\bullet} =\quad \cdots\lra F^i\xra{\
\delta^i\ } F^{i-1}\lra\cdots\lra F^0\lra0\lra\cdots
   \]
where each $\delta^i$ is a homomorphism of degree $0$ and satisfies
$\delta^i(F^i)\subseteq JF^{i-1}$, see \cite[15 and \S4]{Ei}.
Totaling the complex $F\vphantom{I}^{\scriptscriptstyle\bullet}$
one gets a DG $A$-module $F$ with $F^\natural=\coprod_{i=0}^\infty
\shift^iF^i$ and $\dd(F)\subseteq JF$. It comes with a quasi-isomorphism
$F\to N$ and a semi-free filtration $(F^n)_{n\in\BZ}$ with
$F^n{}^\natural=\coprod_{i=0}^n \shift^iF^i$.  If $N$ has a finite
free resolution $G\vphantom{I}^{\scriptscriptstyle\bullet}$, then
$F\vphantom{I}^{\scriptscriptstyle\bullet}$ is isomorphic to as a direct
summand of $G\vphantom{I}^{\scriptscriptstyle\bullet}$ by \cite[8]{Ei},
and hence $N$ is perfect by Theorem \ref{plevel:filtrations}.  Conversely,
if $\level AAM$ is finite, then \ref{levels:finiteness2} implies that
the graded $k$-vector space $\hh{k\lotimes AM}$ has finite  rank.
The isomorphisms
  \[
\hh{k\lotimes AM}\cong \hh{k\lotimes AN}\cong\hh{k\otimes_AF}=
\coprod_{i=0}^\infty\shift^i(k\otimes_AF^i)\,.
  \]
yield $\rank_k\hh{k\lotimes AM}=\sum_{i=0}^\infty\rank_AF^i$, so
$F\vphantom{I}^{\scriptscriptstyle\bullet}$ is a finite free resolution.
 \end{proof}

We conclude with examples that show that for non-positive algebras
condition (b) in Theorem~\ref{level:DGperfection2} cannot be relaxed
significantly; compare with Remark~\ref{level:others}.

  \begin{examples}
Let $k$ be a field and $A$ a DG algebra with 
 \[
A^\natural= k[y,z]/(y^2,yz,z^2)
 \quad\text{and}\quad
\dd^A=0\,,
 \]
where $y$ and $z$ are indeterminates over $k$.  The exact sequence
  \[
0\lra\shift^{|y|}ky\oplus\shift^{|z|}kz \lra A\lra k\lra0
  \]
implies that in every minimal free resolution
$F\vphantom{I}^{\scriptscriptstyle\bullet}$ of the graded
$A^\natural$-module $k$ one has $\rank_AF^i=2^i$ for each $i\ge0$;
by Proposition~\ref{level:DGperfection3}, this yields
 \[
    \level AAk = \infty\,.
 \]

\begin{subchunk}
For $|y|=-1=|z|$ one has $A_0=k$ and $A_{-1}\ne0$.  

Let $M$ be the DG module with $M^\natural=A$ and $\dd(a)=ya$.

The graded $A^\natural$-module $M^\natural$ is free of rank $1$, and
one has $\level AAM=\infty$.
 \end{subchunk}

Indeed, the map $\shift^{-1}k\to M$ sending $\sigma^{-1}(1)$ to $z$
is a quasi-isomorphism of DG $A$-modules, so one has 
$\level AAM=\level AAk=\infty$.

  \begin{subchunk}
For $|y|=0$ and $|z|=-2$ one has $A_0\ne k$ and $A_{-1}=0$.  

Let $M$ be the DG module with $M^\natural=A\oplus\shift^{-1}A$
and $\dd(a,\sigma^{-1}(b))\!=(zb,\sigma^{-1}(ya))$.

The graded $A^\natural$-module $M^\natural$ is
free of rank $2$, and one has $\level AAM=\infty$.
  \end{subchunk}

Indeed, the map $(a,\sigma^{-3}(b))\mapsto\big(ya,\sigma^{-1}(zb)\big)$
is a quasi-isomorphism of DG $A$-modules $k\oplus\shift^{-3}k\to M$,
which yields $\level AAM=\level AA{k\oplus\shift^{-3}k}=\infty$.
  \end{examples}

    \section{A New Intersection Theorem for DG algebras}
    \label{A New Intersection Theorem for DG algebras}

The New Intersection Theorem is a central result in the homological theory
of commutative noetherian rings.  Here we generalize it to certain DG
modules, using Hochster's notion of super height of an ideal $I$ in a
commutative ring $R$:
    \begin{gather*}
     \supheight I = \sup\left\{\height(IS)\left|
      \begin{gathered}
    \text{$R\to S$ is a homomorphism}\\ \text{of rings and $S$ is
    noetherian}
      \end{gathered}
\right\}\right.
 \end{gather*}
Evidently, when $R$ is noetherian one has that $\supheight I\geq
\height I$.

  \begin{theorem}
\label{plevel:intersection}
Let $A$ be a DG algebra with zero differential, let $M$ be a DG module
over $A$, and let $I$ denote the annihilator of $\bigoplus_{n\in\BZ}\HH
nM$ in the ring $A^\flat=\bigoplus_{n\in\BZ} A_n$.

If  $A^\flat$ is commutative and noetherian and is an algebra over a
field, then one has
  \[
\level AAM\geq \supheight I+1\,.
  \]
When the ring $A^\flat$ is Cohen-Macaulay, or when its dimension is at
most $3$, one has
  \[
\level AAM\geq \height I+1\,.
  \]
 \end{theorem}

\begin{Addendum}
\label{plevel:general}
Whenever $A^\flat$ is commutative and noetherian, one has 
\[
\level AAM\geq \supheight I\,.
\]
\end{Addendum}

For the next remark, recall that a ring is just a DG algebra concentrated
in degree $0$, and a DG module over a ring is nothing but a complex.

 \begin{remark}
We recall the statement of the New Intersection Theorem: 

In a bounded complex $P$ of finite free modules with non-zero homology
of finite length over a local ring $R$ one has $P_n\ne0$ for at least
$(\dim R+1)$ values of $n$.

For algebras over a field it is due to Peskine and Szpiro\cite{PS},
P.\ Roberts\cite{Rb1}, and Hochster\cite{Ho2}; this case follows from
Theorem~\ref{plevel:intersection}, because Example~\ref{plevel:card}
yields
   \[
\card\{n\in\BZ\mid P_n\ne0\}\ge\level RRP\,.
   \]
However, we do not recover the New Intersection Theorem over arbitrary
local rings, proved by Roberts; see \cite{Rb1}.  The reason is that in the
proof of Theorem~\ref{plevel:intersection} we need a result from \cite
{dm}, which uses Hochster's \cite{Ho2} big Cohen-Macaulay modules.
 \end{remark}

Over the rings it covers, Theorem~\ref{plevel:intersection} may provide
a significantly sharper bound on heights than the one given by the New
Intersection Theorem, as the difference $\card\{n\var F_n\ne 0\}-\level
AAF$ can be arbitrarily large, even when $R$ is very nice:

  \begin{example}
Let $(R,\fm,k)$ be a regular local ring of dimension $d\ge2$ and $l$
an integer, $l\ge d$.  Let $I$ be an $\fm$-primary ideal minimally
generated by $l$ elements; one always exists.  The Koszul complex $F$
on a minimal generating set for $I$ has
  \[
\level RRF = d + 1
  \quad\text{and}\quad
\card\{n\mid F_n\ne 0\}= l+1\,.
  \]

Indeed, the equality on the right reflects the construction of $F$.
Theorem~\ref{plevel:intersection} yields $\level RRF\geq d+1$,
because $\Ann_R\hh K=I$ and $\height I=d$.  On the other hand,
$F$ is quasi-isomorphic to the complex $D$ obtained by totaling
its Cartan-Eilenberg resolution $C$.  Since $\gldim R = d$, by
\cite[XVII.1.4]{CE} one can choose $C$ to have $d+1$ columns; this yields
a semi-free filtration of $D$ of class at most $d$.
  \end{example}

A consequence of the New Intersection Theorem is that a local
ring with a non-zero module of finite length and finite projective
dimension is Cohen-Macaulay.  Using this and the 
argument above, one can prove a stronger statement:

  \begin{remark}
If $M$ is a complex over a local ring $R$ and the module
$\bigoplus_{n}\HH nM$ has non-zero finite length and finite
projective dimension, then $\level RRM=\dim R+1$.
  \end{remark}

The next result provides an upper bound for $A$-levels that, together
with Theorem \ref{plevel:intersection}, completes the proof of Theorem
\ref{intro:plevel} from the introduction.  

Recall that a graded algebra is said to be \emph{coherent} if finite graded submodules of finitely presented graded modules are also finitely presented.

  \begin{theorem}
\label{plevel:pdim}
Let $A$ be a DG algebra with zero differential, which is left coherent
as a graded ring, and let $M$ be a DG $A$-module.

If $\hh M$ is finitely presented over $A$, then the following
inequalities hold:
    \[
\level AAM\le\pd_A\hh M+1\le\gldim A+1\,.
    \]
   \end{theorem}

  \begin{remark}
The second inequality in Theorem~\ref{plevel:pdim} holds by definition.  

When $A$ is a ring and $M$ is an $A$-module, Krause and
Kussin \cite[2.4]{KK} prove that the first inequality
in Theorem~\ref{plevel:pdim} becomes an equality.
  \end{remark}

In view of the Syzygy Theorem, Theorems \ref{plevel:intersection} and
\ref{plevel:pdim} yield:

    \begin{corollary}
  \label{plevel:poly}
  \pushQED{\qed}%
Let $S$ be a graded polynomial algebra in $c$ indeterminates over a field $k$ with $\dd^{S}=0$.
Each DG $S$-module $M$ with $\rank_k\hh M\ne0,\infty$ satisfies
\begin{equation*}
\level SSM=c+1\,.
\qedhere
    \end{equation*}
   \end{corollary}

  \begin{proof}[Proof of Theorem \emph{\ref{plevel:pdim}}]
It suffices to deal with the case $\pd_A\hh M=p<\infty$.

When $p=0$ the graded $A$-module $\hh M$ is projective.  As $A$ has zero
differential, the cycles of $M$ form a graded $A$-submodule $\ZZ(M)$,
and the canonical surjection $\ZZ(M)\to\hh M$ is an $A$-linear map.
Choosing an $A$-linear splitting $\sigma\col\hh M \to\ZZ(M)$ and composing
it with the inclusion $\ZZ(M)\subseteq M$ one gets a quasi-isomorphism
$\hh M\to M$.  The DG module $\hh M$ is a direct summand of some finite
free graded $A$-module, so in $\dcat A$ one has $\hh M\in\sadd A$.
Thus, one has
   \[
\level AAM=\level AA{\hh M}\le1\,.
   \]

Let now $p$ be a positive integer, and assume that the desired inequality
holds for all DG $A$-modules whose homology has projective dimension
strictly smaller than $p$. In $M$, pick cycles $z_1,\dots,z_s$
whose homology classes generate $\hh M$, set $L=\bigoplus_{i=1}^s
\shift^{|z_i|}A$.  The map of graded $A$-modules $\lambda\col L\to M$
that sends $1\in\shift^{|z_i|}A$ to $z_i$ is a morphism of DG $A$-modules.
In $\dcat A$ it fits into an exact triangle $L\xra{\lambda} M\xra{\mu}
N\to\shift L$.  In the induced exact sequence of graded $A$-modules
    \[
\shift^{-1}\hh M\xra{\ \shift^{-1}\hh{\mu}\ }\shift^{-1}\hh N\lra\hh L
\xra{\hh{\lambda}\ }\hh M\xra{\hh{\mu}\ } \hh N
   \]
the map $\hh{\lambda}$ is surjective by construction.  This implies
$\hh{\mu}=0=\shift^{-1}\hh{\mu}$, so the exact sequence shows that
$\hh N$ is finitely presented over $A$, and one has $\pd_A\hh N=p-1$.
{F}rom the induction hypothesis we now obtain $\level AAN\le p$, hence
we get $\level AAM\le p+1$ from Lemma~\ref{levels:properties}(2).
   \end{proof}

In order to obtain lower bounds in Theorem~\ref{plevel:intersection} on 
$A$-levels of DG $A$-modules we use our recent results on invariants 
of a related structure, which we define next.

  \begin{bfchunk}{Differential modules.}
\label{dg:differential}
Let $R$ be an associative ring.

A \emph{differential $R$-module} is a pair $(D,\delta)$, where $D$ is an
$R$-module and $\delta\col D\to D$ an $R$-linear map with $\delta^2=0$;
the module $\hh D=\Ker(\delta)/\Image(\delta)$ is the \emph{homology}
of $D$.  Every module $M$ supports a differential module $(M,0)$, also
called $M$.

A \emph{projective flag} in a differential $R$-module $D$ is a family
$\Filt Dn_{n\in\BZ}$ of differential $R$-submodules of $D$, such that
for each $n\in\BZ$ the following hold:
    \begin{gather*}
\filt Dn \subseteq \filt D{n+1}, \quad \filt D{-1}=0, \quad
\bigcup_{n\in\BZ}\filt Dn=D\,,\quad\text{and}
   \\
\filt Dn/\filt D{n-1}\cong(P_n,0)
\quad \text{for some projective $R$-module } P_n\,.
    \end{gather*}
In \cite[2.8]{dm} the \emph{projective class} of $D$ is defined to be
the number
   \[
\prclass RD = \inf\left\{l\in\BZ \,\left|\,
   \begin{gathered}
D \text{ admits a projective flag}
   \\
\Filt Dn_{n\in\BZ} \text{ with }\filt D{l}=D
    \end{gathered}
\right\}\right..  \]
   \end{bfchunk}

The statement below is one of the main results of \cite{dm}:

  \begin{bfchunk}{Class Inequality.}
\label{intersection:global}
If $R$ is a commutative noetherian ring, $F$ a finitely
generated differential $R$-module, $D$ a retract of $F$,
and $I=\Ann_R \hh D$, then one has
   \[
\prclass RF \geq \supheight I -1\,,
   \]
with strict inequality when $R$ is an algebra over a field; see \cite[4.2]{dm}. Moreover,
when $\dim R\leq 3$ or $R$ is Cohen-Macaulay \cite[4.1]{dm} yields an inequality
   \[
\prclass RF \geq \height I\,.
   \]
   \end{bfchunk}

There is an obvious parallel between the notion of semi-free filtration
for DG modules and that of projective flag for differential modules.
Under additional conditions we turn it into a direct
comparison, by using the following construction.

  \begin{bfchunk}{DG algebras with zero differential.}
   \label{levels:flattening}
Let $A$ be a DG algebra with $\dd^A=0$.

Using the product of the graded algebra $A$ one defines an associative ring
    \[
\dfm{A} = \bigoplus_{n\in \BZ}A_n\,,
    \]
and then to each DG $A$-module $M$ one associates a differential
$\dfm{A}$-module
    \[
\dfm{M}=\bigg(\bigoplus_{n\in\BZ}M_n\,,\, \bigoplus_{n\in\BZ} (-1)^n
\dd^M_n\bigg)\,.
    \]

The simplicity of the construction notwithstanding, a couple of caveats
may be in order:  The signs appearing in the formula for the differential
of $\dfm{M}$ are necessary to ensure that it is $\dfm A$-linear.
The action of $\dfm A$ on $\dfm M$ respects the obvious internal gradings
of these objects, but the differential of $\dfm M$ need not.
    \end{bfchunk}

  \begin{lemma}
    \label{levels:class}
Let $A$ be a DG algebra with trivial differential.

The assignment $M\mapsto\dfm{M}$ defines an exact functor from the abelian
category of DG $A$-modules to that of differential $\dfm A$-modules.
It preserves finite generation and transforms semi-free filtrations into
projective flags.

If $F$ is a semi-free DG $A$-module of class at most $c$, then one has
   \[
c\ge\prclass{\dfm A}{(\dfm F)}\,.
   \]
  \end{lemma}

  \begin{proof}
The first two assertions are evident from the definitions and
constructions preceding the lemma.  If $\Filt Fn$ is a semi-free
filtration of $F$ with $\filt F{c}=F$, then $((\dfm F){}^n)$ is a
projective flag with $(\dfm F){}^{c}=\dfm F$: this gives the inequality
above.
    \end{proof}

  \begin{proof}[Proof of Theorem \emph{\ref{plevel:intersection}}]
Set $\level AAM=l$. We may assume that $l$ is finite.

Theorem~\ref{plevel:filtrations} shows that $M$ is isomorphic, in
$\dcat A$, to a retract of some finite semi-free DG $A$-module $F$ of
class $l-1$.  Levels do not change under isomorphisms, so we may assume
that $M$ itself is a direct summand of $F$.  In that case $\dfm M$
is a retract of $\dfm F$, and the $\dfm A$-module $\dfm F$ is finite.
Thus, Lemma~\ref{levels:class} yields
   \[
l-1\ge\prclass{\dfm A}{(\dfm F)}\,.
   \]
To finish the proof, invoke the Class Inequality from
\ref{intersection:global}.
    \end{proof}

\section{Levels and semi-simplicity}
\label{Levels and semi-simplicity}

In this section we analyze levels of DG modules related to two classical
notions of length for modules over rings.

  \begin{bfchunk}{Restricted lengths of modules.}
 \label{lol:definition}
Let $S$ be a ring and $\cc$ a finite non-empty set of simple
$S$-modules.  A \emph{$\cc$-filtration} of an $S$-module $N$ is a
sequence of submodules
\[
0=\filt{N}{0} \subsetneq \filt{N}{1}\subsetneq \cdots
\subsetneq \filt{N}{l}=N
\]
where every $\filt{N}{i}/\filt{N}{i-1}$ is isomorphic to a direct sum of
modules from $\cc$.

The \emph{$\cc$-length of $N$}, denoted $\slength S{\cc}N$, is the
largest integer $l$ for which $N$ has a $\cc$-filtration with $N^l=N$;
when no such filtration exists we set $\slength S{\cc}N=\infty$.

The \emph{$\cc$-Loewy length of $N$}, denoted $\slol S{\cc}N$, is the
least integer $l$ for which $N$ has a $\cc$-filtration with $N^l=N$;
when none exists we set $\slol S{\cc}N=\infty$.

  \begin{subchunk}
\label{lol:additive}
For each exact sequence $0\to N'\to N\to N''\to 0$ of $S$-modules one
has
\begin{align}
   \label{length:additive}
    \tag{\ref{lol:additive}.1}
\slength S{\cc}N
&=\slength S{\cc}{N'} + \slength S{\cc}{N''}
\\
   \label{lol:subadditive}
     \tag{\ref{lol:additive}.2}
\slol S{\cc}N
&\leq \slol S{\cc}{N'} + \slol S{\cc}{N''}
   \\
\slol S{\cc}{N}
&\ge\max\{\slol S{\cc}{N'},\slol S{\cc}{N''}\}
   \label{lol:split}
     \tag{\ref{lol:additive}.3}
\intertext{with equality in \eqref{lol:split} when the sequence splits.
{F}rom here one gets}
     \label{lol:S}
     \tag{\ref{lol:additive}.4}
     \slol S{\cc}N &\leq \slol S{\cc}S\,.
    \end{align}
      \end{subchunk}

  \begin{sublemma}
    \label{lol:length}
Let $N$ be an $S$-module.

If $\slength S{\cc}N$ is finite, then so is $\slol S{\cc}N$.

When $N$ is noetherian the converse holds as well.
    \end{sublemma}

  \begin{proof}
By definition, one has $\slol S{\cc}N\le \slength S{\cc}N$, whence
the first assertion. For the converse, assume $\slol S{\cc}N=l<\infty$
and fix a filtration
    \[
0=\filt{N}{0} \subsetneq \filt{N}{1}\subsetneq \cdots
\subsetneq \filt{N}{l} = N
    \]
    where each $\filt{N}{i}/\filt{N}{i-1}$ a direct sum of simple  
modules from $\cc$.  When
    $N$ is noetherian each direct sum is finite, so $N$ has finite $ 
\cc$-length.
  \end{proof}

Recall that $S$ is \emph{semi-local} if it has finitely many isomorphism
classes of simple modules; equivalently, if $S/\fn$, where $\fn$ is the
Jacobson radical, is semi-simple.

  \begin{subchunk}
    \label{lol:semilocal}
Let $S$ be a semi-local ring, let $\fn$ denote its Jacobson radical, and
let $\cc$ contain representatives of every isomorphism class of simple
$S$-modules.  The $\cc$-length and $\cc$-Loewy length of $N$
are then equal to their classical counterparts, denoted $\length_SN$
and $\lol SN$, respectively.  Furthermore, one has
    \[
\lol SN = \inf\{n\in\BN \mid \fn^n N =0\}
               =\inf\{n\in\BN \mid (0:\fn^n)_N =N\}\,.
    \]
    \end{subchunk}
      \end{bfchunk}

Now we return to DG algebras.  As in \ref{augmentation:non-negative}
and \ref{augmentation:non-positive}, we use canonical augmentations
$\eps\col A\to\HH 0A$ to identify $\HH0A$-modules with DG $A$-modules.

Parts (3) and (4) of the next result contain Proposition \ref{intro:length}
from the introduction.

  \begin{theorem}
\label{klevel:estimates}
Let $A$ be a non-negative DG algebra.  Set $S=\HH 0A$, let $\cc$ be a
finite set of simple $S$-modules, and set $k=\bigoplus_{N \in\cc}N$.

For each DG $A$-module $M$ the following statements hold.
\begin{enumerate}[{\quad\rm(1)}]
\item
$\level A{k}M = \level A{\cc}M$.
\item
$\level{A}{k}M\geq
\max_{n\in \BZ}\{\,\slol {S}{\cc}{\HH nM}\,\}$.
\item
The numbers $\level A{k}M$ and $\slength {S}{\cc}
{\hh M}$ are finite simultaneously; when they are, they
are linked by the following inequality:
   \[
\level{A}{k}M \leq \sum_{n\in\BZ} \slol {S}{\cc}
{\HH nM}\,.
   \]
\item
If $N$ is an $S$-module and $\slength{S}{\cc}N$ is finite, then
one has
   \[
\level{A}{k}N=\slol{S}{\cc}N\,.
   \]
  \end{enumerate}
    \end{theorem}

  \begin{proof}
(1) This follows from the equality $\smd{\sadd k} = \smd{\sadd{\cc}}$
in $\dcat S$.

(2) We may assume $\level A{k}M=l<\infty$.  By (1) and
\ref{layers:symmetry}, this means that $M$ is a direct summand of
a complex $L\in{\sadd{\cc}}^{\star l}$.  When $l=1$ one may assume
that each $L_i$ is a finite direct sum of modules from $\cc$; then so
is each $\HH iL$, and thus $\slol S{\cc}{\HH iL}\le1$.  For $l\ge 2$
there is an exact triangle $L'\to L\to L''\to \shift L'$ in $\dcat S$
with $\level A{k}{L'}\le l-1$ and $\level A{k}{L''}\le1$.  It induces
an exact sequence $\HH i{L'}\to \HH iL\to\HH i{L''}$ of $S$-modules.
Induction and sub-additivity, see \eqref {lol:subadditive}, yield $\slol
S{\cc}{\HH iL}\le l$.  Thus, one gets $\slol S{\cc}{\HH iM}\le l$.

(3) It is easy to verify that the subcategory of DG $A$-modules $M$,
such that $\slength S{k}{\hh M}$ is finite, is thick.  It contains $k$,
hence also $\thick{A}{k}$.  Thus, when $\level A{k}M$ is finite, so is
$\slength S{\cc}{\hh M}$, see \ref{levels:inheritance}.  One now has
    \[
     \level A{k}M\leq \sum_{n\in\BZ} \level S{k}{\HH nM}
    \]
from Proposition~\ref{truncations:non-negative}. To finish, for each
$S$-module $N$ of finite $\cc$-length we prove:
  \[
\level SkN \leq \slol S{C}N\,.
  \]
Indeed, $N$ has a filtration by submodules with $\filt{N}{0}=0$, $
\filt{N}{l} = N$ for $l=\slol S{C}N$, and each $\filt{N}{i}/\filt{N}{i-1}$
a direct sum of modules from $\cc$.  The sums are finite because $\slength
S{\cc}N$ is, so each $\filt{N}{i}/\filt{N}{i-1}$ has $\cc$-level $1$. Thus
\ref{levels:finiteness1} yields $\level{S}{k}N \leq l$.

(4) This is a formal consequence of (2) and (3).
    \end{proof}

We pause to remark that the finiteness hypothesis in Theorem~\ref{klevel:estimates}(4), as well as the
noetherian hypothesis in Lemma~\ref{lol:length} are essential:

  \begin{example}
When $k$ is a simple $S$-module the module $N=k^{(\BN)}$ satisfies
   \[
\slol SkN=1<\infty=\level SkN=\slength SkN\,.
   \]
    \end{example}

More precise conclusions can be made when $\HH0A$ is semi-simple.

  \begin{theorem}
\label{klevel:weight}
Let $A$ be a DG algebra satisfying one of the conditions below:
\begin{enumerate}[{\quad\rm(a)}]
 \item
$A$ is non-negative.
 \item
$A$ is non-positive and $\dd(A_0)=0$.
  \end{enumerate}
Set $S=\HH 0A$ and assume that the ring $S$ is semi-simple.

For a DG $A$-module $M$ the number $\level A{S}M$ is finite if and
only if the $S$-module $\hh M$ is finitely generated; when it is, the
following inequality holds:
  \[
\level{A}{S}M \leq \card\{n\in\BZ \mid \HH nM\ne 0\}\,.
  \] \end{theorem}

  \begin{proof}
The semi-simple ring $S$ is noetherian.  Let $\cc$ be a set of
representatives of the isomorphism classes of simple $S$-modules, and
note that $\cc$ is finite.

If $\level ASM$ is finite, then the $S$-module $\hh M$ is noetherian;
see \ref{levels:finiteness2}.

Conversely, assume that the graded $S$-module $\hh M$ is finitely
generated.  It is isomorphic to a direct sum of suspensions of modules
from $\cc$, so for each $n\in\BZ$ one has $\HH nM\neq0$ if and only
if $\level SS{\HH nM} =1$.  The desired inequality follows from
Proposition~\ref{truncations:non-negative}(2) in case (a) and from Proposition
\ref{truncations:non-positive}(2) in case (b).
    \end{proof}

The next two results link the structure of a semi-local ring to
$k$-levels.

  \begin{proposition}
  \label{klevel:artinian}
Let $S$ be a left artinian ring and $\fn$ its Jacobson radical.

If $M$ is a complex of $S$-modules with $\hh M$ finitely generated,
then one has
  \[
\level S{S/\fn}M\le\lol SS\,.
  \]
  \end{proposition}

  \begin{proof}
Replacing $M$,  if necessary, by a quasi-isomorphic complex, we may assume that $M$ is a bounded complex of finite $S$-modules.  Since $\fn^l=0$ for $l=\lol SS$, setting $\filt{M}i=\fn^{l-i}M$ we get a filtration of
$M$ by subcomplexes with  $\filt{M}{0} = 0$, $\filt{M} {l} = M$ and $\filt{M}{i}/\filt{M}{i-1}$ a bounded complex of finite semi-simple $S$-modules for $i=0,\dots,l$.  In $\dcat S$ such a complex is {quasi-isomorphic} to its homology.  The latter, being a bounded complex of finite semi-simple modules with trivial differential, has $k$-level $1$ by definition. Thus \ref{levels:finiteness1} now yields $\level S{S/\fn}M\le
l$, as desired.
    \end{proof}

  \begin{proposition} \label{klevel:one}
For a local ring $(R,\fm,k)$ the following conditions are equivalent.
  \begin{enumerate}[\quad\rm(i)]
\item
  The ring $R$ is regular.
\item
  There exists a finite free complex of $k$-level $1$.
\item
  The Koszul complex $K$ on a minimal set of generators of $\fm$ has
  $k$-level $1$.
  \end{enumerate}
   \end{proposition}

  \begin{proof}
When $R$ is regular in $\dcat R$ one has $K\simeq k$, hence
$\level{R}kK=1$.  

Let $F$ be a finite free complex with $\level{R}kF=1$.  By definition,
in ${\dcat R}$ one then has $F\simeq V$, where $V$ is a complex of
$k$-vector spaces with $\rank_k\hh V$ finite non-zero.  As $k$ is a retract
of $\shift^sV$ for some $s\in\BZ$, in $\dcat R$ it is also a retract
of $\shift^sF$, and so is quasi-isomorphic to a finite free complex.
This means that the $R$-module $k$ has finite projective dimension,
hence $R$ is regular; see \cite[2.2.7]{BH}.
    \end{proof}

We finish the section with examples of strict inequalities in the
preceding results.

  \begin{example}
Let $(R,\fm,k)$ be a local ring and $K$ the
Koszul complex on a minimal set of generators of $\fm$.
Each module $\HH nK$ is a vector space over $k$, and so
has $k$-level $1$; see Theorem~\ref{klevel:estimates}(4).
When $R$ is singular, one gets the inequality below
\[
\level{R}{k}K>1=\max_{ n\in\BZ}\{\slol Rk{\HH nK}\}\,,
\]
so the inequality in Theorem~\ref{klevel:estimates}(2) can
be strict.  On the other hand, $\HH nK\ne0$ holds precisely
when $n$ satisfies $0\le n\le \edim R-\depth R$, so one has
\[
\sum_{n\in\BZ} \slol{R}{k}{\HH nK}=
\card\{n\in\BZ \mid \HH nK\ne 0\}=\edim R-\depth R+1\,.
\]
The number on the right hand side is independent of $\lol RR$,
so by fixing one and varying the other it is easy to conjure
artinian local rings for which a strict inequality holds in
Theorem~\ref{klevel:estimates}(3), Theorem \ref{klevel:weight}, or
Proposition~\ref{klevel:artinian}.
  \end{example}

\section{Perfect DG modules over exterior algebras}
\label{Perfect DG modules over exterior algebras}

Our goal is to prove a slightly enhanced version of Theorem
\ref{intro:bgg} in the introduction. It is an algebraic analogue of,
and can be used to deduce, \cite[4.4.5]{AP}, which is a statement
about the toral rank of spaces.

  \begin{theorem}
\label{bgg:weight}
Let $k$ be a field and $\Lambda$ a DG $k$-algebra with $\dd^{\Lambda}=0$
and $\Lambda^\natural$ an exterior algebra on $c$ alternating
indeterminates of positive odd degrees.

For every perfect DG $\Lambda$-module $N$ with $\hh N\ne 0$ one then has
\[
\level{\Lambda}kN=c+1\leq\card\{n\in\BZ\mid\HH nN\ne0\}\,.
\]
  \end{theorem}

It will be deduced from Corollary \ref{plevel:poly} by using equivalences
between certain thick subcategories of $\dcat\Lambda$ and $\dcat S$,
where $S$ is a graded polynomial ring in $c$ indeterminates over $k$.
The prototype of such results is a classical theorem of J.~Bernstein,
I.~M.~Gelfand, and S.~I.~Gelfand~\cite{BGG}, dealing with subcategories
of the derived category of graded modules over these graded algebras.
The situation here is different.  We provide a self-contained treatment,
as none of the results that we have located in the literature
covers it with the detail and in the generality that we need; see
Remarks~\ref{bgg:history} and \ref{bgg:DG}.

  \begin{chunk}
    \label{bgg:data}
Let $k$ be a field and let $c$ be a non-negative integer.

\smallskip

Let  $\Lambda$ be the DG algebra with $\dd^\Lambda=0$ and
$\Lambda^\natural$ an exterior algebra on alternating indeterminates
$\xi_1,\dots,\xi_c$ of positive odd degrees.

\smallskip

Let $S$ denote the DG algebra with $\dd^S=0$ and $S^\natural$ a polynomial
ring on commuting indeterminates $x_1,\dots,x_c$, with $|x_i|=-|\xi_i|-1$
for $i=1,\dots,c$.

\smallskip

Set $(-)^*=\Hom k-k$, viewed as a functor on the category of complexes
of $k$-vector spaces; see Section \ref{Levels of DG modules}.  The DG
algebras $S$ and $\Lambda$ are graded commutative, so for every DG module
$L$ over either one of them the complex $L^*$ of $k$-vector spaces carries
a canonical structure of DG module over the same DG algebra.  \end{chunk}

  \begin{chunk}
 \label{bgg:X}
As $S$ and $\Lambda$ are graded commutative with elements of odd
degree squaring to zero, the DG algebra $\Lambda\otimes_kS$ has the
same property.  In particular,
  \[
\delta=\sum_{h=1}^c\xi_h\otimes x_h\in(\Lambda\otimes_kS)_{-1}
\quad\text{satisfies}\quad \delta^2=0\,.
  \]

Let $E$ be a DG $(\Lambda\otimes_kS)$-module.  An elementary calculation
shows that the map
    \[
\dd\col E\to E \quad\text{defined by}\quad \dd(e)=\dd^E(e)+\delta\cdot e
    \]
for all $e\in E$ is $k$-linear of degree $-1$, and satisfies $\dd^2=0$
and $\dd(ae)=(-1)^{|a|}a\dd(e)$ for all $a\in\Lambda\otimes_kS$.  Thus,
$E^\delta=(E^\natural,\dd)$ is a DG $(\Lambda\otimes_kS)$-module.

\smallskip

For the DG $\Lambda$-module $\Lambda^*$, see \ref{bgg:data}, and the
DG $S$-module $S$, the tensor product $\Lambda^*\otimes_kS$ is naturally
a DG module with zero differential over $\Lambda \otimes_kS$.  Set
    \[
X=(\Lambda^*\otimes_kS)^\delta\,.
    \]
\end{chunk}

The next result provides the last ingredient needed in the proof of Theorem \ref{bgg:weight}. The DG module structures of $\Hom{\Lambda}{X^{*}}N$ over $S$ and of $X^{*}\lotimes SM$ over $\Lambda$ are those induced by $X^{*}$, as explained in \ref{delta-calculus}.

 \begin{theorem}
  \label{bgg:theorem}
Let $X$ be the DG module described above. 

The functor $\Hom{\Lambda}{X^{*}}-$ induces an exact functor $\SH\col \dcat{\Lambda}\to \dcat{S}$. 

The functor $\ST=X^{*}\lotimes S-\col \dcat S\to \dcat{\Lambda}$ is left adjoint to $\SH$. 

This pair of functors restricts to inverse equivalences of triangulated categories:
    \begin{equation}
  \label{ }
   \begin{gathered}
\xymatrixcolsep{2pc} \xymatrix{
\dcat{\Lambda}\ar@{<-}[r]<.6ex>^-{\ST}
\ar@{}[d]|-{\vgu}
&\dcat{S}\ar@{<-}[l]<.6ex>^-{\SH}
\ar@{}[d]|-{\vgu}
\\
\thick{\Lambda}{k}\ar@{<->}[r]^{\equiv}
\ar@{}[d]|-{\vgu}
&\thick{S}{S}
\ar@{}[d]|-{\vgu}
\\
\thick{\Lambda}{\Lambda}\ar@{<->}[r]^{\equiv}
&\thick{S}{k}
}
 \end{gathered}
  \end{equation}
For $d=|\xi_1|+\cdots+|\xi_c|$ there are isomorphisms in $\dcat S$ and $\dcat\Lambda$, respectively:
    \begin{equation}
    \label{bgg:calculations}
\begin{aligned}
\SH(k)&\simeq S
\\
\SH(\Lambda)&\simeq \shift^d k 
  \end{aligned}
\qquad\text{and}\qquad
\begin{aligned}
\ST(S)&\simeq k
\\
\ST(k)&\simeq\shift^{-d}\Lambda 
  \end{aligned}
    \end{equation}
    \end{theorem}
    
The theorem is proved at the end of this section.  The thick
subcategories that appear in its statement admit very explicit
descriptions.

  \begin{remark}
    \label{bgg:finiteness}
For each DG $S$-module $M$ the following hold:
\begin{align}
\label{bgg:SS-finite}
M\in\thick{S}{S} &\iff  M\simeq M' \text{ with $M'{}^\natural$ finite  
free over } S
\\
\label{bgg:SSS-finite}
\tag{\ref{bgg:finiteness}.1$'$}
                  &\iff \hh M \text{ is finite over }S\,.
\\
\label{bgg:Sk-finite}
M\in\thick{S}{k} &\iff \hh M \text{ is finite over }k\,.
\\
\intertext{\indent For each DG $\Lambda$-module $N$ the following hold:}
\label{bgg:LL-finite}
N\in\thick{\Lambda}{\Lambda} &\iff N\simeq N' \text{ with $N'{}^ 
\natural$ finite free over }\Lambda\,.
\\
\label{bgg:Lk-finite}
N\in\thick{\Lambda}{k} &\iff \hh N \text{ is finite over }k
\\
\label{bgg:LLL-finite}
\tag{\ref{bgg:finiteness}.4$'$}
                  &\iff \hh N \text{ is finite over }\Lambda\,.
  \end{align}

Indeed, \eqref{bgg:SS-finite} and \eqref{bgg:LL-finite} are special cases
of Theorem~\ref{level:DGperfection2}, while \eqref{bgg:Sk-finite} and
\eqref{bgg:Lk-finite} are special cases of Theorem~\ref{klevel:weight}.
When $\level SSM$ is finite the graded $S$-module $\hh M$ is noetherian
by \ref{levels:finiteness2}: this establishes one direction of
\eqref{bgg:SSS-finite}.  Conversely, when $\hh M$ is finite from
Theorem~\ref{plevel:pdim} and the Syzygy Theorem one gets $\level
SSM\le \gldim S=c+1$.  Finally, \eqref{bgg:LLL-finite} is evident as
$\rank_k\Lambda$ is finite.
    \end{remark}

 \begin{bfchunk}{The Koszul DG module.}
For use in the proof of Theorem~\ref{bgg:theorem}, we collect homological
properties of the DG module $X$ introduced in \ref{bgg:X}.

  \begin{subchunk}
 \label{Xexact}
Let $\eta^{\Lambda}\col k\to\Lambda$ denote the structure map and
$\eps^S \col S\to k$ the canonical augmentation.  The following map
of DG $S$-modules is a quasi-isomorphism:
  \[
\pi\col X=(\Lambda^*\otimes_kS)^\delta\xra{\ (\eta^{\Lambda}){}^*
\otimes\eps^S\ } k\otimes_kk=k\,.
  \]
Indeed, bigrading the complex of $k$-vector spaces underlying $X$
by assigning to the indeterminates $\xi_i$ and $x_i$ homological degrees
$1$ and $0$, respectively, one obtains the graded Koszul complex on
$x_1,\dots,x_c\in S$; its homology is equal to $k$.
  \end{subchunk}

  \begin{subchunk}
 \label{Xcoexact}
Let $\eta^S\col k\to S$ denote the structure map and $\eps^{\Lambda}\col\Lambda\to k$ the
canonical augmentation.  The following map is a quasi-isomorphism of DG $\Lambda$-modules:
\[
\iota\col k=k\otimes_kk\xra{\ (\eps^{\Lambda}){}^*\otimes\eta^S\ }
(\Lambda^*\otimes S)^\delta=X\,.
  \]
As $\pi\iota=\idmap^k$ and $\pi$ is a quasi-isomorphism, $\iota$ is 
one as well.  The induced map 
  \[
\rho=\Hom{k}{\iota}{k}\col X^*\lra k
  \]
is a quasi-isomorphism of DG $\Lambda$-modules.
  \end{subchunk}

 \begin{subchunk}
\label{delta-calculus}
For every DG $\Lambda$-module $N$ the complex $\Hom{\Lambda}{X^*}N$ of $k$-vector spaces has a canonical  structure of DG $S$-module, isomorphic to $\Hom{k}{S^{*}}{N}^{\delta}$. 

Multiplication with $\delta$ annihilates $\Hom{\Lambda}{\Lambda\otimes_kS^*}k{}^\natural$.  This observation shows that the following natural $S$-linear isomorphisms are morphisms of DG $S$-modules:
\[
\Hom{\Lambda}{X^*}{k} \cong\Hom{\Lambda}{(\Lambda\otimes_kS^*)^{-\delta}}k
=\Hom{\Lambda}{\Lambda\otimes_kS^*}k \cong S\,.
\]

For each DG $S$-module $M$ the complex $X^{*}\otimes_SM$ of $k$-vector spaces has a  structure of DG $\Lambda$-module, isomorphic to $(\Lambda^{*}\otimes_{k}M)^{\delta}$. Thus, $X^{*}\otimes_{S}-$ is a functor from DG $S$-modules to DG $\Lambda$-modules; it induces a functor $X^{*}\lotimes S-\col \dcat S\to \dcat \Lambda$.
\end{subchunk}
 
\begin{subchunk}
  \label{XoverS}
The following functors preserve quasi-isomorphisms of DG $S$-modules
 \[
\Hom SX-\,,\quad
(X\otimes_S-)\,,
\quad\text{and}\quad
\big(-\otimes_S\Hom{\Lambda}{X^*}{\Lambda}\big)\,,
 \]
Indeed, by Theorem \ref{level:DGperfection2} and
\ref{semi-free} it suffices to note that $X^\natural$ and
$\Hom{\Lambda}{X^*}{\Lambda}^\natural$ are finite free over $S$.
For the first module this is clear; for the second one has
  \[
(\Hom{\Lambda}{X^*}{\Lambda})^\natural\cong\Hom k{S^*}{\Lambda}
\cong\Lambda\otimes_kS\,.
  \] \end{subchunk}

\begin{subchunk}
\label{XoverL}
The following functors preserve quasi-isomorphisms of DG
$\Lambda$-modules:
 \[
\Hom {\Lambda}{X^*}-
 \quad\text{and}\quad
\Hom {\Lambda}-X\,.
 \]

Indeed, since $X^*{}^\natural\cong \Lambda\otimes_kS^*$ as graded
$\Lambda$-modules,  $X^*$ is a direct summand of a semi-free DG
$\Lambda$-module by Proposition \ref{level:semi-projective}.  Now apply
\ref{semi-free} to $X^{*}$ and $X\cong X^{**}$.
 \end{subchunk}

\begin{subchunk}
\label{exterior:self-dual}
There exists an isomorphism of DG $\Lambda$-modules:
$\Lambda^*\cong\shift^{-d}\Lambda$; in particular, the functor $\Hom
{\Lambda}-{\Lambda}$ preserves quasi-isomorphisms of DG $\Lambda$-modules.

Indeed, $(\Lambda^*)_{-d}$ contains a unique $k$-linear map
that sends $\xi_1\cdots\xi_c$ to $1_k$.  It defines an element
$\omega\in(\shift^{d}(\Lambda^*))_0$.  The homomorphism 
$\Lambda\to\shift^{d}(\Lambda^*)$ of graded $\Lambda$-modules, given by
$\lambda\mapsto \lambda\omega$, is easily seen to be injective, and
hence is bijective.
 \end{subchunk}
  \end{bfchunk}

  \begin{proof}[Proof of Theorem~\emph{\ref{bgg:theorem}}]
  The functor $\Hom{\Lambda}{X^*}-$ preserves quasi-isomorphisms, by \ref{XoverL}, so it induces an exact functor $\SH\col\dcat\Lambda\to\dcat S$. It is clear that the functor $\ST=X^{*}\lotimes S-$ is a left adjoint to $\SH$. Let   
\[
\sigma\col \mathsf{Id}^{\dcat S} \to \SH\ST
\qquad\text{and}\qquad
\lambda\col \ST\SH\to\mathsf{Id}^{\dcat \Lambda}\,,
    \]
respectively, be the unit and the counit of the adjunction. Thus, for each $M\in\dcat S$ and $N\in\dcat \Lambda$ there are morphisms of DG modules over $S$ and $\Lambda$, respectively:
\[
M\xra{\sigma(M)} \Hom{\Lambda}{X^{*}}{X^{*}\lotimes SM}
\qquad\text{and}\qquad
X^{*}\lotimes S \Hom{\Lambda}{X^{*}}N \xra{\lambda(N)} N\,.
\]
To prove that $\SH$ and $\ST$ restrict to inverse equivalences between  $\thick{\Lambda}{k}$ and $\thick{S}{S}$ it suffices to show $\lambda(k)$  and $\sigma(S)$  are quasi-isomorphisms.

Consider the chain of morphisms of DG $\Lambda$-modules
\[
 X^{*}\lotimes S \Hom{\Lambda}{X^{*}}k \cong X^{*}\lotimes SS\cong X^{*}\xra{\rho} k\,,
\]
where the first isomorphism comes from \ref{delta-calculus}. By the same token, the DG $S$-module $\Hom{\Lambda}{X^{*}}k$
is semi-free, so in the formula above one may replace $\lotimes S$ with $\otimes_{S}$. Now a direct verification shows that the composite map equals $\lambda(k)$. Since $\rho$ is a quasi-isomorphism, see \ref{Xcoexact}, so is $\lambda(k)$.

In the following chain of morphisms of DG $S$-modules the last one is from \ref{delta-calculus}:
\[
S\xra{\sigma(S)}\Hom{\Lambda}{X^{*}}{X^{*}}\xra{\Hom{\Lambda}{X^{*}}{\rho}}\Hom{\Lambda}{X^{*}}k \cong S.
\]
A direct computation shows that the composite map is $\idmap^{S}$. Since $\rho$ is a quasi-isomorphism so is $\Hom{\Lambda}{X^{*}}\rho$, by \ref{XoverL}; hence $\sigma(S)$ is a quasi-isomorphism.

We have proved that $\SH$ and  $\ST$  restrict to inverse equivalences between $\thick {\Lambda}k$  and $\thick SS$. We now verify the isomorphisms in \eqref{bgg:calculations}.

We obtain an isomorphism $\SH(k)\simeq S$ from \ref{delta-calculus}. Using the quasi-isomorphism $\rho$ of DG $\Lambda$-modules from \ref{Xcoexact} and the exactness of the functor $\Hom{\Lambda}-{\Lambda}$, see
\ref{exterior:self-dual}, one gets quasi-isomorphisms
 \[
\SH(\Lambda)=\Hom{\Lambda}{X^*}{\Lambda}\simeq
\Hom{\Lambda}k{\Lambda}\cong\shift^dk
 \]
of complexes of vector spaces.  It yields $\SH(\Lambda)\simeq\shift^dk$ in $\dcat S$,
see Proposition \ref{truncations:non-positive}(1).

Since $\lambda\col \ST\SH\to \idmap^{\dcat \Lambda}$ is an isomorphism on $\thick{\Lambda}k$, it follows from the calculations above that there are quasi-isomorphisms of DG $\Lambda$-modules
\[
\ST(S)\simeq \ST\SH(k)\simeq k \quad\text{and}\quad \ST(k)\simeq \shift^{-d} \ST\SH(\Lambda) \simeq \shift^{-d}\Lambda
\]
  
Finally, since $\SH$ restricts to an equivalence $\thick{\Lambda}k\equiv \thick SS$ and $\SH(\Lambda)\simeq \shift^{d}k$ holds, further restriction yields an equivalence $\thick{\Lambda}{\Lambda}\equiv\thick Sk$.
\end{proof}

  \begin{proof}[Proof of Theorem \emph{\ref{bgg:weight}}]
Recall the hypothesis: $N$ is a DG $\Lambda$-module in
$\thick\Lambda\Lambda$, with $\hh N\ne 0$.  The equivalences of categories
$\SH$ from Theorem~\ref{bgg:theorem} imply that $\rank_k\hh{\SH(N)}$
is  finite and non-zero, and yields the first equality below:
    \[
\level{\Lambda}kN = \level SS{\SH(N)}=c+1\,.
    \]
The second one is given by Corollary~\ref{plevel:poly}.  Finally,
from Theorem~\ref{klevel:weight}(a) we get
 \[
\card\{n\in\BZ\mid \HH nN\ne 0\} \geq \level{\Lambda}kN\,.
 \qedhere
 \]
    \end{proof}

  \begin{remark}
   \label{bgg:history}
Several publications deal with equivalences of subcategories of
the derived categories of DG modules $\dcat S$ and $\dcat\Lambda$:
see \cite[II.7]{Ca:ln}, \cite[App.]{AM} and the bibliography of
the latter.  None of these is applicable to the present situation, for
reasons having to do with specific choices of gradings, restrictions on
the characteristic of $k$, focus on different subcategories, or reliance
on tools from analysis.
   \end{remark}

  \begin{remark}
\label{bgg:DG}
Techniques in \cite{DGI1} allow for a different approach to
Theorem~\ref{bgg:theorem}.

The first step is to interpolate the endomorphism DG algebra $E=\Hom {\Lambda}XX$ between
$\Lambda$ and $S$.  Using \cite[4.10]{DGI1} one can then prove that the functors
 \[
(-\lotimes{\ce}k) \col \dcat E\to \dcat \Lambda \quad\text{and}\quad
\Rhom{\Lambda}k-\col \dcat\Lambda\to \dcat E
 \]
 induce an equivalence of triangulated categories between $\dcat\Lambda$ and a certain
 subcategory of $\dcat E$.  It is well-known that the homothety map $S\to E$ of DG
 algebras is a quasi-isomorphism. It remains to track this subcategory of $\dcat E$ under
 the equivalence $\dcat S\equiv\dcat E$.  The necessary calculations are similar to
 those used above to prove Theorem~\ref{bgg:theorem}, so we have presented a direct
 approach.
  \end{remark}

\section{The conormal rank of a local ring}
\label{The conormal rank of a local ring}

For each local ring we introduce a numerical invariant that corresponds to the maximal
rank of a free direct summand of the conormal module of a graded algebra. 
We start by describing a language convenient for the discussion.

  \begin{chunk}
\label{cohen}
Let $(R,\fm,k)$ be a local ring.

A \emph{local presentation} of $R$ is simply an isomorphism of rings
$R\cong Q/I$, where $(Q,\fq,k)$ is a local ring. We say that such a
presentation is \emph{minimal} if $\edim R=\edim Q$, that is, if the
ideal $I$ is contained in $\fq^2$.

A local presentation $R\cong Q/I$ is \emph{regular} if the ring
$Q$ is regular.  If it is not minimal, then choosing an element
$x\in I\smallsetminus\fq^2$ one obtains a regular presentation
$R\cong(Q/Qx)/(I/Qx)$ with $\edim(Q/Qx)= \edim Q-1$.  Iterating this
procedure, one sees that every regular presentation can be factored
through a minimal one.

Regular presentations exist when $R$ is essentially of finite type
over a field, or when $R$ is complete: the latter case comes from
Cohen's Structure Theorem, which provides a presentation
$R\cong Q/I$ with a complete regular local ring $Q$.
  \end{chunk}

Recall that for an $R$-module $M$ the number $\frank_RM$, called
the  \emph{free rank} of $M$ over $R$, is defined to be the maximal
rank of a free direct  summand of $M$.

  \begin{chunk}
\label{conormal:rank}
We define the \emph{conormal free rank} of $R$ to be the number
\[
\cnr R = \sup\left\{\frank_{\wh R}(I/I^2)
       \left|
\begin{gathered}
   \wh R\cong Q/I\text{ is a minimal}
\\
  \text{regular presentation}
  \end{gathered}\right.\right\}\,.
\]
  \end{chunk}

Some estimates on the new invariant are easy to come by:

  \begin{lemma}
  \label{conormal:freesummands}
Let $R\cong Q'/I'$ be a minimal presentation.
  \begin{enumerate}[\quad\rm(1)]
   \item
The following equality holds:  $\cnr R=\cnr{\wh R}$.
  \item
The following inequality holds: $\cnr R \ge \frank_{R}(I'/I'{}^2)$.
  \item
If $I'=Q'\bsx+I''$ and $\bsx$ is a $(Q'/I'')$-regular sequence of length
$r$, then
   \[
\cnr R \ge r+\frank _{\wh R}(\ov I/\ov I{}^2)\,,
  \]
where $\ov I$ denotes the ideal $I'/(Q'\bsx)$ in $\ov Q=Q'/(Q'\bsx)$.
   \end{enumerate}
    \end{lemma}

  \begin{proof}
(1)  The desired equality follows directly from the definition.

(2) We may assume $R$ and $Q'$ are complete, as the induced isomorphism $\wh R\cong {\wh Q}'/{\wh I}'$ is a minimal presentation and one has $\frank_{\wh R}({\wh I}'/{\wh I}'{}^2)\ge \frank_R(I'/I'{}^2)$. Choose a minimal regular presentation $Q'\cong Q/J$.  The composition $Q\to Q'\to \wh R$ induces a minimal regular presentation $R\cong Q/I$ and a surjective homomorphism of $R$-modules $I/I^2\to I'/I{}'^2$, which yields
    \[
\frank_R(I/I^2)\ge\frank_R(I'/I'{}^2)\,.
    \]

(3) This follows from a well known fact recorded in
\ref{conormal:regularsequence} below.
  \end{proof}

  \begin{chunk}
  \label{conormal:regularsequence}
Let $R\cong Q'/I'$, where $I'=Q'\bsx+I''$ and $\bsx=\{x_1,\dots,x_r\}$.
Set $\ov Q=Q'/(Q'\bsx)$ and $\ov I=I'/(Q'\bsx)$.  If $\bsx$ is
$Q'/I''$-regular, then the sequence of $R$-modules
  \[
0\lra R^r\xra{\ \varkappa\ }I'/I'{}^2\lra \ov I/\ov I{}^2\lra 0
  \]
where $\varkappa$ maps the $i$th element of the standard basis of
$R^r$ to $x_i+{I'}^2$ is split exact.

Indeed, for $r=1$ this is  in \cite[Cor. pp 458]{Ho1}; the general case
follows by iteration.
  \end{chunk}

One says that $R$ is \emph{equicharacteristic} if
$\operatorname{char}(k)\cdot R=0$. We show that for such rings there
exists a reasonable notion of \emph{conormal module}.  Whether this is
so in general is related to Grothendieck's Lifting Problem, discussed
by Hochster in \cite{Ho1}.

  \begin{proposition}
\label{conormal:unique}
Let $(R,\fm,k)$ be an equicharacteristic local ring.

If $\wh R\cong Q/I$ and $\wh R\cong Q'/I'$ are minimal regular presentations with $Q,Q'$ equicharacteristic,  the $\wh R$-modules $I/I^2$ and $I'/I'^2$ are isomorphic.  In particular, one has
  \[
\cnr R = \frank_{R}(I/I^2)\,.
  \]
  \end{proposition}

  \begin{proof}
Passing to completions, we may assume they are complete.  The fiber product $Q\times_{R}Q'$ then is a
complete equicharacteristic local ring.  Choosing a  minimal regular
presentation of  $Q\times_{R}Q'$, see \ref{cohen}, we get a
commutative diagram of surjective homomorphisms of complete
equicharacteristic  local rings
  \[
\xymatrixrowsep{1.5pc}
\xymatrixcolsep{1.8pc}
\xymatrix{
   &&Q\ar@{->}[dr]&
\\
P  \ar@{->}[r]
  \ar@/^1.8pc/[urr]^-{\pi}
 \ar@/_1.8pc/[drr]_-{\pi'}
  &Q\times_{R}Q'\ar@{->}[ur]\ar@{->}[dr]
  & & {\  R\   }
\\
  &&Q'\ar@{->}[ur]&}
  \]

As $P$ and $Q$ are regular, the ideal $\Ker(\pi)$ is generated
by elements $x_1,\dots,x_r$ that form part of a regular system of
parameters for $P$.  One has
  \[
r=\dim P-\dim Q=\dim P-\edim Q=\dim P-\edim R\,.
  \]

By Cohen's Structure Theorem, the complete equicharacteristic local
ring $P$ is a ring of formal power series over $k$ in $e=\dim P$ formal
variables.  It follows that one may find elements $x_{r+1},\dots,x_e$
in $P$ such that $P=k[\![x_1,\dots,x_e]\!]$; consequently,
$Q=k[\![x_{r+1},\dots,x_e]\!]$. There is thus an ideal $J''\subseteq
(x_{r+1},\dots,x_e)$ in $P$ such that $J''Q=I$ and $J=(x_1,\dots,x_r)+J''$ where
$J=\Ker(P\to R)$.  Now \ref{conormal:regularsequence} yields
  \[
J/J^2\cong I/I^2\oplus{R}^r
  \]
as $R$-modules. By symmetry, a similar isomorphism holds with $I'$ in
place of $I$, so we get an isomorphism of finitely generated modules
over the complete ring $R$:
  \[
I/I^2\oplus{R}^r\cong I'/I'^2\oplus{R}^r\,.
  \]

{F}rom the Krull-Remak-Schmidt Theorem one gets $I/I^2\cong I'/I'^2$.
  \end{proof}

\section{DG algebra models for Koszul complexes}
\label{DG algebra models for Koszul complexes}

In this section $(R,\fm,k)$ denotes a local ring and $K$ the Koszul
complex on a minimal set of generators of $\fm$.  We prove two results
that play a significant role in the proof of our main theorem, given
in the next section.

  \begin{theorem}
\label{koszul:models}
Let $A$ be a non-negative DG algebra, linked to $K$ by a  sequence
of quasi-isomorphisms of  DG algebras, and let
    \[
\SJ\col \dcat K\xra{\ \equiv\ }\dcat A
    \]
be the induced equivalence of categories.  In $\dcat A$ there is then an isomorphism
    \[
\SJ(K\otimes_Rk) \simeq \bigoplus_{n\ges0}\shift^n{\HH0A}^{\binom en}
\quad\text{where}\quad e=\edim R\,.
    \]
  \end{theorem}

A special case is worth noting: in $\dcat K$ one has $K\otimes_Rk \simeq \bigoplus_{n\ges 0}\shift^n k^{\binom en}$.
  
The second result involves the conormal free rank of $R$, an invariant of the ring denoted
$\cnr R$ which is introduced and discussed in Section~\ref{The conormal rank of a
local ring}.

  \begin{theorem}
\label{koszul:structure}
Let $(R,\fm,k)$ be a local ring, and set
\[
e=\edim R\,,\quad d=\depth R\,,\quad c=\cnr R\,.
\]
Let $\Lambda$ denote the exterior algebra $\bigwedge_k(\shift k^c)$
with $\dd(\Lambda)=0$.

There exist quasi-isomorphisms of DG algebras linking $K$ and $\Lambda\otimes_kB$, where $B$ is a DG algebra with $B_0=k$, $\rank_k B<\infty$, $\dd(B_1)=0=\dd(B_{2})$, and
    \[
\sup\{i\in\BN\var\HH iB\ne0\}=e-d-c\,.
    \]
  \end{theorem}

  To prove Theorems~\ref{koszul:models} and \ref{koszul:structure} we use DG algebra resolutions of a special type.

  \begin{chunk}
 \label{models}
Let $(Q,\fq,k)$ be a local ring. A \emph{semi-free $\varGamma$-extension}
of $Q$ is a DG algebra $Q\la X\ra$, where $X$ is a set of divided powers
indeterminates with $X_n$ finite for each $n\ge1$ and empty for $n \le0 $.
Thus, the graded algebra underlying $Q\la X\ra$ is a tensor product of
the exterior algebra on the free $Q$-module with ordered basis $X_ \odd$
and the divided powers algebra on the free $Q$-module with ordered basis
$X_\even $.

For $y\in X_\even$ we let $y^{(n)}$ denote the $n$th divided power of $y$.
A $Q $-basis of $Q\la X\ra$ is given by all products $x_{i_1}\cdots
x_{i_r}y_{j_1}^{(n_1)} \cdots y_{j_s}^{(n_s)}$ with $x_{i_u}\in X_\odd$
and $i_1<\cdots<i_r$, with $y_{i_v}\in X_\even$ and $j_1<\cdots<j_s$,
and with $r,s,n_v\ge0$.

We let $X^{(\ges 2)}$ denote the set of monomials with
$r+n_1+\cdots+n_s\ge2$.
    \end{chunk}

  \begin{chunk}
    \label{models:quadratic}
Let $R\cong Q/I$ be a regular presentation; see \ref{cohen}.  An
\emph{acyclic closure} of $R$ over $Q$ is a quasi-isomorphism
$\vf\col Q\la X\ra\to R$ of DG algebras satisfying
    \[
\dd(X) \subseteq I + \fq X + QX^{(\ges 2)}
    \]
Every presentation can be extended to an acyclic closure, see
\cite[\S6.3]{Av:ifr}.
    \end{chunk}

  \begin{lemma}
    \label{koszul:klifting}
Let $\pi\col R\to k$ denote the natural surjection, $\rho\col
R\to\wh R$ the completion map, and $(Q,\fq)\to \wh R$ a minimal regular
presentation.  Let $E$ be the Koszul complex on a minimal generating set
of $\fq$ and $\eps\col E\to k$ the canonical augmentation.

If $\vf\col Q\la X\ra\to\wh R$ is an acyclic closure, then there is a
morphism $\vk\col Q\la X\ra\to E$ of DG $\varGamma$-algebras making the
following diagram commute.
    \[
\xymatrixrowsep{2pc} \xymatrixcolsep{2pc} \xymatrix{
R\ar@{->}[r]^-{\rho}\ar@{->}[d]_-{\pi} &{\wh R}\ar@{->}[d]_-{\wh\pi}
\ar@{<-}[r]_\simeq^-{\vf} & Q\la X\ra\ar@{->}[d]^{\vk}
    \\
k\ar@{=}[r] &k\ar@{<-}[r]_-{\eps}^-{\simeq}
    & {E} } \]
Furthermore, each such morphism satisfies $\vk(Q\la X\ra)_{\ges
1}\subseteq \fq E$.
    \end{lemma}

  \begin{proof}
The map $\eps$ is a quasi-isomorphism because $Q$ is regular.  The
existence of $\vk$ follows from general principles: $Q\la X\ra$ is free as
an algebra with divided powers over $Q$, the $Q$-algebra $E$ has a system
of divided powers, and $E\to k$ is a quasi-isomorphism.  To finish the
proof it suffices to show that $\vk(Q\la X_{\les n}\ra)_{\ges 1}\subseteq
\fq E$ holds for each integer $n\geq 0$.  We proceed by induction on $n$.

The desired inclusion holds trivially for $n=0$ because one has $Q\la
X_{\les 0}\ra=Q$ by definition.  Assuming that the inclusion holds for
some $n\geq 0$, one has
    \begin{align*}
\dd\vk(X_{n+1})
                &=\vk\dd(X_{n+1}) \\
                &\subseteq \vk(IQ_0 + \fq X_n + (X_{\les n})^{(\ges 2)}) \\
                &\subseteq \fq^2E_0+\fq \vk(X_n)+(\vk(X_{\les n}))^{(\ges 2)}\\
                &\subseteq \fq^2 E
    \end{align*}
Indeed, the first equality holds because $\vk$ is a morphism of
complexes. The first inclusion comes from the formula in \ref{models:quadratic}.
The second one is a consequence of the hypothesis that the presentation
is minimal and the fact that $\vk$ is a homomorphism of algebras with
divided powers.  The last inclusion results from the induction hypothesis
and the identity $(ay)^{(n)}=a^ny^{(n)}$ for $a\in Q$.

As an algebra with divided powers, $Q\la X_{\les n+1}\ra$ is generated
over $Q\la X_{\les n}\ra$ by the set $X_{\les n+1}$.  Thus, the
inclusion established above implies $\dd(\vk(Q\la X_{\les n+1}\ra)_{\ges
1})\subseteq \fq^2 E$.  A result of Serre \cite[Ch.I\ V, App.\ I, Prop.\
3]{Se}, see \cite[4.1.6(2)]{Av:ifr}, now yields $\vk(Q\la X_{\les
n+1}\ra)_{\ges 1}\subseteq \fq E$, and the induction step is complete.
    \end{proof}

  \begin{proof}[Proof of Theorem \emph{\ref{koszul:models}}]
Recall that $K$ denotes the Koszul complex on a minimal set of generators
of $\fm$.  In view of Proposition~\ref{dg:h0actions} and \ref{morphisms:objects}, to prove the theorem
it suffices to construct a specific non-negative DG algebra $A$ with
$\HH0A=k$ and to produce a chain of quasi-isomorphisms linking $K $ and
$A$, such that for the induced equivalence of categories $\SJ\col\dcat
K\to \dcat A$ there is an isomorphism
  \begin{equation}
 \label{iso}
\SJ(K\otimes_Rk)\cong \bigoplus_{n\ges0}\shift^nk^{\binom en}
  \end{equation}
Let $\rho\col R\to \wh R$ be the completion map and $Q\to \wh R$ a minimal regular 
presentation. Set $A=k\otimes_QQ\la X\ra$, where $\vf\col Q\la X\ra\to\wh R$ is an
acyclic closure, see \ref{models:quadratic}. Set  $\wh K=K\otimes_R\wh R$, let
$E$ be the Koszul complex on a minimal generating set of the maximal
ideal $\fq$ of $Q$, and let $\eta\col k\to k \otimes_QE$ denote the
structure map.  We claim that the following diagram of morphisms of DG
$Q$-algebras commutes:
    \[
\xymatrixrowsep{2pc} \xymatrixcolsep{.75pc} \xymatrix{
k \ar@{=}[rr]\ar@{<-}[d]
    && k \ar@{=}[r]\ar@{<-}[d]
     & k \ar@{=}[rrr]\ar@{<-}[d]
      &&& k \ar@{=}[rrr]\ar@{<-}[d]
       &&& k \ar@{<-}[d]\ar@/^2pc/[dd]^-{\eta}
    \\
K_{\vphantom{Q}}{\vphantom{\wh K}}\ar@{->}[rr]^-{K\otimes_R\rho}_-
{\simeq} \ar@{->}[d]_-{K\otimes_R{\pi}}
    &&{\wh K}_{\vphantom Q} \ar@{->}[d]_-{\wh K\otimes_{\wh R}{\wh\pi}}
\ar@{=}[r]
     & E\otimes_Q\wh R\vphantom{\wh K}\ar@{->}[d]_-{E\otimes_Q{\wh
\pi}}\ar@{<-}[rrr]_-\simeq^-{E\otimes_Q\vf}
      &&& E\otimes_Q Q\la X\ra{\vphantom{\wh K}}\ar@{->}[d]_{E\otimes_Q
\vk}\ar@{->}[rrr]_-\simeq^-{\eps\otimes_QQ\la X\ra}
       &&& A\vphantom{\wh K}\ar@{->}[d]_{k\otimes_Q\vk}
    \\
K\otimes_R k{\vphantom{\wh K}}\ar@{->}[rr]^-{\simeq}
    &&{\wh K}\otimes_{\wh R}k \ar@{=}[r]
     & E\otimes_Qk{\vphantom{\wh K}}\ar@{<-}[rrr]^-\simeq_-{E\otimes_Q
\eps}
      &&& E\otimes_Q E{\vphantom{\wh K}}\ar@{->}[rrr]^-\simeq_-{\eps
\otimes_QE}
       &&&k\otimes_QE{\vphantom{\wh K}}
     }
    \]
Indeed, Lemma~\ref{koszul:klifting} shows that the squares commute
and yields $\vk(Q\la X\ra_{\ges 1})\subseteq \fq E$; this implies
$(k\otimes_Q\vk)(A_{\ges 1})=0$, so the triangle on the right commutes
as well.

In the bottom left square one has  quasi-isomorphisms because $\rho$
is flat and the $R$-modules $\hh K$ and $\hh{K\otimes_Rk}$ have finite
length.  The other quasi-isomorphisms hold because $E$ and $Q\la X\ra$
are bounded below  complexes of free $Q$-modules.

The middle row in the diagram is a sequence of quasi-isomorphism of DG
algebras linking $K$ with $A$.  Let $\SJ\col \dcat K\to \dcat A$ be the
equivalence of triangulated categories induced by it.  In view of the
commutativity of the bottom part of the diagram, repeated application
of \ref{morphisms:objects} gives in $\dcat A$ an isomorphism
    \[
\SJ(K\otimes_Rk) \simeq k\otimes_QE
    \]
The triangle containing $\eta$ implies $A_{\ges1}\cdot (k\otimes_QE) =0$,
while the inclusion $\dd(E)\subseteq\fq E$ yields $\dd(k\otimes_QE)=0$.
Thus, one has an isomorphism of DG $A$-modules
    \[
k\otimes_QE \cong  \bigoplus_{n\ges0}  \shift^nk^{\binom en}
    \]

Concatenating the isomorphisms above one obtains \eqref{iso}, as desired.
  \end{proof}

The next proof is similar to that of \cite[2.1]{Iy}; it uses an idea
from \cite{An}.

  \begin{proof}[Proof of Theorem \emph{\ref{koszul:structure}}]
Let $\wh R\cong Q/I$ be a minimal regular presentation such that one
has $\frank_{\wh R}(I/I^{2})=c$. Choose elements $\bsa= a_{1},\dots,a_{c}$
in $I$ whose images form a basis for a free summand of $I/I^{2}$. The
Koszul complex on the elements $\bsa$ is an exterior algebra $Q\la
x_{1},\dots,x_{c}\ra$ with the $\bsx$ indeterminates in degree $1$ and
differential defined by $\dd(x_{i})=a_{i}$. By Nakayama's lemma $\bsa$
can be extended to a minimal generating set for the ideal $I$, so this
Koszul complex can be extended to acyclic closure $Q\la X\ra\xra{\simeq}
\wh R$ with $\{x_{1},\dots,x_{c}\}$ a subset of $X_{1}$.

We use the notation of Lemma \ref{koszul:klifting}. It follows
from the proof of  Theorem~\ref{koszul:models} that there exists a
quasi-isomorphism of DG algebras
    \[
K\to \wh K =  E\otimes_Q\wh R\xla{\ E\otimes_Q\vf\ } E\otimes_QQ\la X\ra
    \]
We claim that one has an isomorphism of DG algebras
    \[
E\otimes_QQ\la X\ra\cong Q\la z_1,\dots,z_c\ra \otimes_Q E\la Y\ra\,,
    \]
where $|z_{i}|=1$ and $\dd(z_i)=0$ for $i=1,\dots,c$, and
$Y_n=\varnothing$ for $n\le0$.

Indeed, there are derivations $\theta_{1},\dots,\theta_{c}$ of the
DG $\varGamma$-algebra $Q\la X\ra$ such that $\theta_{i}(x_{1})=1$
and $\theta_{i}(x)=0$ for all  $x\in X_{1}\setminus \{x_{1}\}$;
see \cite[6.2.7]{Av:ifr} and also \cite[1.4]{Iy}. Replacing
$\theta_{i}$ by  $\theta_{i} - x_{i}\theta_{i}^{2}$ one may assume
in addition $\theta_{i}^{2}=0$ holds. These derivations extend
to $E$-linear derivations $\wt\theta_{1},\dots,\wt\theta_{c}$
of the DG $\varGamma$-algebra $E\la X\ra=E\otimes_{Q}Q\la X\ra$
with $(\wt\theta_{i})^{2}=0$. Since $I \subseteq \dd(E_{1})$
holds there are elements $e_{1},\dots,e_{c}$ in $E_{1}$ such that
$\dd(e_{i})=a_{i}$ for each $i$. Evidently, the elements $z_{i}=
x_{i}-e_{i}$ are cycles of degree $1$ in $E\la X\ra$ satisfying
$\wt\theta_{i}(z_{j})=\delta_{ij}$. Set $Y'=\{x - z_{1}\wt\theta_{1}(x)\var
x\in X\setminus \{x_{1}\}\}$. These are indeterminates over $E$ and the
induce a bijective morphism of DG $\varGamma$-algebras
  \[
Q\la z_{1}\ra  \otimes_{Q} E\la Y'\ra  \to E\la X\ra\,.
  \] 
It is not hard to verify that the derivations
$\wt\theta_{2},\dots,\wt\theta_{c}$ restrict to derivations on $E\la
Y'\ra$. Iteratiing the procedure one gets the desired isomorphism of
DG algebras.

Since $E\la Y\ra$ is semi-free as a DG $E$-module, the quasi-isomorphism
$\eps\col E\to k$ induces a quasi-isomorphism of DG algebras
    \[
E\la Y\ra=E\la Y\ra\otimes_EE\xra{E\la Y\ra\otimes_E\eps} E\la
Y\ra\otimes_Ek=C
    \]
where $C=k\la Y\ra$.  In turn, it induces the quasi-isomorphism below:
    \[
Q\la z_1,\dots,z_c\ra\otimes_QE\la Y\ra \xra{\simeq} Q\la
z_1,\dots,z_c\ra\otimes_QC \cong \Lambda\otimes_kC
    \]

The description of the acyclic closure $Q\la X\ra$, see \ref{models:quadratic}, shows that $C$ is non-negative, $\rank_kC_n$ is finite for each $n$, $C_0=k$, and
$\dd(C_{1})=0=\dd(C_{2})$ holds. Set $b=e-d-c$ and note the equalities
    \[
\sup\{i\in\BN\var\HH iC\ne0\}= \sup\{i\in\BN\var\HH iK\ne0\}-c=b
    \]
The first one comes from the isomorphisms
$\hh{K}\cong\hh{\Lambda\otimes_kC}\cong\Lambda\otimes_k\hh{C}$.
The second one is from the Auslander-Buchsbaum Formula.  It is easy to
check that \[ D=\quad \cdots \to C_{b+2} \to C_{b+1}\to \dd(C_{b+1})\to
0 \] is a DG ideal of $C$ with $\hh D=0$.  Thus, $B=C/D$ is a DG
$k$-algebra with $B_n=0$ for $n>b$ and the canonical surjection $C\to B$
is a quasi-isomorphism of DG algebras.  It induces a quasi-isomorphism
of DG algebras $\Lambda\otimes_kC\to\Lambda\otimes_kB$.

Now we have produced a DG algebra $B$ that has the required properties,
and we have linked $K$ with $\Lambda\otimes_kB$ by a chain of
quasi-isomorphisms of DG algebras.
  \end{proof}

\section{Loewy length of homology of perfect complexes}
\label{Loewy length of homology of perfect complexes}

In this section we prove the following result; it contains Theorem~\ref {intro:loewy}.

  \begin{theorem}
\label{conormal:main}
Let $(R,\fm,k)$ be a local ring. Each finite free complex $F$ of  $R$-modules with
$\hh F\ne 0$ satisfies the following inequalities:
    \[
\sum_{n\in\BZ}\lol R{\HH nF} \geq  \level{R}kF \geq \cnr R + 1\,.
    \]
  \end{theorem}

Recall that the \emph{closed fiber} of a local homomorphism $(P,\fp)\to Q$ is the local ring
$Q/\fp Q$. The following corollary contains Theorem~\ref{intro:fibres} from the
introduction.

  \begin{corollary}
\label{loewy:fibres}
If $\wh R$ is the closed fiber of a flat local homomorphism $P\to Q$, then
    \[
\sum_{n\in\BZ}\lol R{\HH nF} \geq s+1\,,
    \]
where $s=\edim P-\edim Q+\edim R$.  In particular, one has $\fm^s\ne0$.
  \end{corollary}

  \begin{proof}
Set $p=\edim P$, $q=\edim Q$, and $r=\edim R$. Let $\bsx=x_1,\dots,x_p$ be
a generating set for the maximal ideal $\fp$ of $P$; thus, $\wh R=Q/\bsx
Q$. Reindexing, if necessary, we may assume $x_i\in \fq \setminus \fq^2$
for $1\leq i \leq q - r$ and $x_i\in \fq^2$ otherwise.

Set $Q'=Q/(x_1,\cdots,x_{q-r})$ and $I'=(x_{q-r+1},\dots,x_p)Q'$. One gets a minimal presentation $\wh R\cong Q'/I'$ with $I'/I'^2\cong(\fp/\fp'')\otimes_kR$, where $\fp''= (x_{1},\dots,x_{q-r})P+\fp^{2}$.  Thus
the $\wh R$-module $I'/I'^2$ is free of rank $p-q+r=s$, so the theorem
applies.
   \end{proof}

Theorem~\ref{conormal:main} has an analogue for complexes of finite injective dimension.

  \begin{corollary}
\label{conormal:injective}
Let $C$ be a bounded complex of injective $R$-modules with $\hh C$
finitely generated and non-zero. The following inequalities then hold:
    \[
\sum_{n\in\BZ}\lol R{\HH nC}\ge\level{R}kC \geq \cnr R + 1\,.
\]
  \end{corollary}

  \begin{proof}
We may assume $\lol R{\hh C}$ is finite. Since the $R$-module $\hh C$ is finitely generated,
this assumption implies that $\length_R\hh C$ is finite, see \ref{lol:length}, so
Theorem~\ref{klevel:estimates}(3)  yields the desired upper bound on $\level{R}kC$.

Let $E$ be the injective hull of $k$. Since $\length_R{\hh C}$ is finite, one may replace $C$ with a quasi-isomorphic complex if necessary, and assume that it is a finite complex with each $C_i$ a finite direct sum of copies of $E$. Thus the complex $F=\Hom RCE$ is finite free over $\wh R$ and the functor $\Hom R-E\col{\dcat R}^\op\to \dcat{\wh R}$ is exact, so
    \[
\level RkC \geq \level {\wh R}kF
    \]
holds by Lemma~\ref{levels:properties}(6).  Combining this with Theorem
\ref{conormal:main}, we get the desired lower bound on $\level{R}kC$.
\end{proof}

Given this corollary, it is clear that one has also an analogue of Corollary~\ref{loewy:fibres} for complexes, and, in particular, for modules, of finite injective dimension.

The preceding results are optimal, in the sense that all the inequalities
involved may become equalities, as the following example demonstrates.

  \begin{example}
\label{loewy:examples}
Let $f_1,\dots,f_p$ be a regular sequence in $k[\![y_1,\dots,y_q] \!]$.
The local ring $R=k[\![y_1,\dots,y_q]\!]/(f_1,\dots,f_p)$ is the closed
fiber of the flat homomorphism
  \[
k[\![x_1,\dots,x_p]\!]\lra k[\![y_1,\dots,y_q]\!]
  \]
of complete $k$-algebras that sends $x_i$ to  $f_i$ for $i=1,\dots p$.

Take $f_i\in\fq^2$ for $i=1,\dots,p$, and $F$ the Koszul complex on the
images in $R$ of $y_1,\dots,y_q$.  One then has $\HH 1F\cong k^p$ and
$\hh F\cong\bigwedge\HH1F$, whence the equalities in the next display,
while Theorem~\ref{conormal:main} and \ref {conormal:rank} give the
inequalities:
  \[
p+1=\sum_{n\in\BZ}\lol R{\HH nF} \geq \level{R}kF \geq \cnr R + 1\ge p+1
  \]
Thus, equalities hold throughout. As one has $p=\edim P$ and $\edim Q=q=\edim R$,
 the inequality in Corollary~\ref{loewy:fibres} is sharp.

Take $p=q$, let $E$ be the injective hull of the $R$-module $k$, and
form the complex $C=\Hom RFE$ of injective $R$-modules.  The isomorphisms
  \[
\HH iC\cong {\Hom R{\HH{-i}F}E}\cong \Hom
k{{\textstyle\bigwedge}^i_k(k^p)}k
  \]
show that all three quantities in the formula of Corollary
\ref{conormal:injective} are equal to $p+1$.
  \end{example}

  \begin{proof}[Proof of Theorem \emph{\ref{conormal:main}}]
Several DG algebras are used in the argument.  Set
 \[
c=\cnr R   
\]
and fix the following notation for the duration of the proof:
    \begin{xalignat*}{1}
K\quad & \text{is the Koszul complex on a minimal generating set for
the ideal }\fm\,.
    \\
\Lambda\quad & \text{is the exterior algebra over $k$ on $c$ variables
of degree $1$, with } \dd^\Lambda=0\,.
    \intertext{By Theorem~\ref{koszul:structure} there exists a DG
$k$-algebra $B$ with $\rank_kB<\infty$, such that} A=&\Lambda\otimes_kB
\text{\ is linked to $K$ by a sequence of quasi-isomorphisms}.
    \end{xalignat*}
The argument hinges on a sequence of exact functors of triangulated
categories
    \begin{equation*}
\xymatrixcolsep{2.2pc}
    \xymatrix{
\dcat R \ar@{->}[r]^{\SK}
     &\dcat K\ar@{->}[r]^-{\SJ}_-\equiv
       &\dcat A\ar@{->}[r]^-{\SII}
        &\dcat \Lambda
}
  \end{equation*}
where the first category is the derived category of $R$-modules, the
rest are derived categories of DG modules over DG algebras.  The functors
are as follows:
    \begin{xalignat*}{1}
\SK&=K\otimes_R-\, \text{; it is exact because the $R$-module $K$
is flat}\,.
    \\
\SJ&\quad \text{is the equivalence induced by the quasi-isomorphisms of
Theorem~\ref{koszul:structure}}.
    \\
\SII&=\iota_{*}\text{ where $\iota\col\Lambda\to \Lambda\otimes_kB$
is the inclusion of DG $k$-algebras}\,.
    \end{xalignat*}
The flow chart below captures the structure of the argument:

The assertion of the theorem is obtained by concatenating the sequence of
(in)equalities displayed in the zigzagging line in the middle.

The column on the left indicates the category in which a given numerical
invariant is computed, and displays functors between such categories.

The column on the right tracks conditions needed for the final equality.
\smallskip\noindent
    \[
\xymatrixcolsep{.1pc}
\xymatrixrowsep{.8pc}
\xymatrix{
    &&&{\sum_{n}}\lol R{\HH nF}\ar@{}[dd]|-{\vge}_{\textstyle(1)\qquad}
      &&&&& \text{$\hh F\ne0$ \& $F$ perfect}\ar@{=>}[dd]^{\qquad\textstyle(1')}
\\
\\
\dcat R \ar@{->}[dd]_{\SK}
     &&& \level{R}{k}{F}\ar@{}[dd]|-{\vge}_{\textstyle(2)\qquad}
        &&&&&  \level{R}{R}{F}\ne 0,\infty \ar@{=>}[dd]^{\qquad
\textstyle(2')}
\\
\\
\dcat K \ar@{->}[dd]_{\SJ}^{\equiv}
     &&&\level{K}{\SK(k)}{\SK(F)}\ar@{=}[dd]_{\textstyle(3)\qquad}
        &&&&& \level{K}{K}{\SK(F)}\ne 0,\infty\ar@{=>}[dd]^{\qquad
\textstyle(3')}
\\
\\
\dcat A \ar@{->}[dd]_{\SII} 
     &&&\level{A}{\SJ\SK(k)}{\SJ\SK(F)} \ar@{=}
         [rr]^-{\textstyle(4)}
                           && \level{A}{k}{\SJ\SK(F)}
\ar@{}[dd]|-{\vge}_{\textstyle(5)\qquad}
        &&& \level{A}{A}{\SJ\SK(F)}\ne 0,\infty\ar@{=>}[dd]^{\qquad\textstyle(5')}
\\
\\
\dcat \Lambda
  &&&&&\level{\Lambda}{k}{\SII\SJ\SK(F)}\ar@{=}[dd]
      &&& \level{\Lambda}{\Lambda}{\SII\SJ\SK(F)}\ne 0,\infty
          \ar@{=>}[dd]^{\qquad\textstyle(6)}
\\
&&&&& {\phantom{xxx}}
\\
&&&&& c+1
      &&{}\ar@{=>}[llu]_-{\textstyle(7)}
       & \text{$\hh{\SII\SJ\SK(F)}\ne0$ \& $\SII\SJ\SK(F)$ perfect}
\\
\\
}
\]
What follows are step by step directions for climbing down the chart.

\smallskip

(1) comes from Theorem~\ref{klevel:estimates}(3) and \ref{lol:semilocal}.

\smallskip

(1$'$) holds by Example \ref{plevel:card}.

\smallskip

(2) holds because $\SK$ is an exact functor; see Lemma~\ref
{levels:properties}(6).

\smallskip

(2$'$) is seen as follows: The exactness of $\SK$ implies, as above,
that $\level KK{\SK(F)}$ is finite; it is non-zero because one has
$\hh{K\otimes_RF}\ne0$, see the end of the proof.

\smallskip

(3) and (3$'$) hold because $\SJ$ is an equivalence of categories;
see Proposition~\ref{levels:basechange2}(1).

\smallskip

(4) holds because one has $\smd{\sadd{\SJ\SK(k)}}=\smd{\sadd{k}}$, by 
Theorem \ref{koszul:models}.

\smallskip

(5) holds because $\SII$ is an exact functor and $\SII(k)=k$; see
Lemma~\ref{levels:properties}(6).

\smallskip

($5'$) is seen as follows: $\level{\Lambda}{\Lambda}{\SII\SJ\SK(F)}$ is non-zero 
since $\hh{\SII\SJ\SK(F)}\cong\hh{\SK(F)}\ne0$; it it finite because
$\level{A}{A}{\SJ\SK(F)}$ is and Corollary~\ref{levels:basechange3}
applies.

\smallskip

(6) holds by definition.

\smallskip

(7) is Theorem \ref{bgg:weight}.

\smallskip

It remains to show that $\hh{K\otimes_RF}$ is not equal to zero.  Since $F$ 
is homologically bounded, the number $i=\inf\{n\in\BZ\var \HH nF\ne0\}$ is finite 
and $F$ is quasi-isomorphic to a complex $L$ with $L_j=0$ for all $j<i$.  
We then get
 \[
\HH i{K\otimes_RF}\cong\HH i{K\otimes_RL}\cong\HH 0K\otimes_R\HH iL
=k\otimes_R\HH iF\ne0
 \]
by using flatness, right exactness of tensor products, and Nakayama's Lemma.
    \end{proof}

\section{Complete intersection local rings}
\label{Complete intersection local rings}

We have thus far focused mainly on complexes admitting finite free resolutions. 
In this final section we turn to complete intersection local rings, over which we 
extend Theorem~\ref{conormal:main} to a statement applying to all homologically 
finite complexes.

\begin{bfchunk}{Complexity.}
Let $(R,\fm,k)$ be a local ring and $M$ a homologically finite complex of 
$R$-modules. A \emph{minimal free resolution} of $M$ is a quasi-isomorphism 
$F\xra{\simeq} M$ where $F$ is a complex of free $R$-modules whose 
differential satisfies $\dd(F)\subseteq \fm F$. Such a resolution exists 
and is unique up to isomorphism of complexes; see \cite{Rb}. 

The \emph{complexity} of $M$ over $R$ is the number
\[
\cxy M = \inf\left\{d\ge 0 \left |
\begin{gathered}
\text{there exists an integer $a>0$ such that} \\
\text{$\rank_{R}(F_{n})\leq a n^{d-1}$ for all $n\gg 0$}
\end{gathered}
\right.
\right\}
\]
Evidently,  $\cxy M=0$ holds if and only if $M$ is perfect.
  \end{bfchunk}

\begin{chunk}
Let $(R,\fm,k)$ be a local ring and $\wh R$ its $\fm$-adic completion. 

The ring $R$ is said to be \emph{complete intersection} if for some 
regular presentation $\wh R\cong P/I$, see \ref{cohen}, the conormal
module $I/I^2$ is free over $\wh R$.  This is equivalent to $I$ being 
generated by a $P$-regular sequence; see \cite[2.2.8]{BH}.  Such
$R$ satisfy
\[
\cnr R =\codim R\,,
\]
where $\cnr R$ is the conormal free rank of $R$; see Lemma~\ref{conormal:freesummands}(3), and $\codim R$ denotes the 
\emph{codimension} of $R$, that is, the difference $\edim R-\dim R$.
\end{chunk}

The following result contains Theorem~\ref{intro:ciall} from the Introduction. 
Over complete intersections it extends Theorem~\ref{conormal:main}, from 
which it is deduced by using a result from \cite{Av:vpd}.  A different approach
to its proof is given in \cite{AI:leeds}.

\begin{theorem}
\label{ci:main}
Let $R$ be a complete intersection local ring and $M$ a complex of 
$R$-modules with $\hh M$ finite and non-zero.  One then has inequalities
\[
\sum_{n\in\BZ}\lol R{\HH nM}\geq \level RkM \geq \codim R - \cxy M + 1\,.
\]
\end{theorem}

The number $\codim R - \cxy M$ in the inequality above is non-negative; 
see \ref{cxy:cibound}.

For the proofs we recall basic facts about complexity.

\begin{chunk}
Let $R$ be a local ring, $M$ a complex of $R$-modules with $\hh M$ finitely 
generated, and let $F$ be a minimal free resolution of $M$.

\begin{subchunk}
\label{cxy:modules}
Fix an integer $s$ such that $\HH nM=0$ for $n\ge s$ and set 
$C=\HH s{F_{\ges s}}$.  Evidently, $\shift^{-s}F_{\ges s}$ is a minimal free 
resolution of $C$, and hence $\cxy C=\cxy M$.
\end{subchunk}

\begin{subchunk}
\label{cxy:cibound}
If there exists a local presentation $\wh R\cong Q/I$ such that $I$ is generated by a regular sequence and  $\wh R\otimes_{R}M$ is perfect over $Q$, then the following inequality holds:
\[
\cxy M \leq \dim Q - \dim R\,.
\]
This follows from \cite[3.2(3)]{Av:vpd}, in view of \ref{cxy:modules}; see also \cite[3.10]{SS}. 
\end{subchunk}
\end{chunk}

When the ring $R$ is a complete intersection, \cite[3.6]{Av:vpd} 
provides a converse to \ref{cxy:cibound} for modules. We extend
that result to complexes.

\begin{proposition}
\label{cxy:deformation}
Let $R$ be a complete intersection with an infinite residue field and $M$ 
a complex of $R$-modules with $\hh M$ finitely generated. 

There exists then a minimal presentation $\wh R\cong Q/J$, such that $J$ 
is generated by a regular sequence of length $\cxy M$ and 
$\wh R\otimes_{R}M$ is perfect over $Q$.
\end{proposition}

\begin{proof}
Set $\cxy M=d$, and let $F$ be a minimal free resolution of $M$. The complex 
$\wh R\otimes_{R}F$ is a minimal free resolution of $\wh R\otimes_{R}M$ over 
$\wh R$, and so $\cxy[\wh R]{(\wh R\otimes M)}=d$ holds. Thus, passing to 
$\wh R$, one may assume that the ring $R$ is complete.

Set $s=\max\{n\mid \HH nM\ne 0\}$ and $C=\HH s{F_{\ges s}}$. Applying 
\cite[3.6]{Av:vpd} to a minimal regular presentation of $R$, see \ref{cohen}, 
one obtains a minimal presentation $R\cong Q/J$ with $J$ generated by a 
$Q$-regular sequence of length $d$ and $C$ perfect over $Q$.

Since $\cxy C=d$, see \ref{cxy:modules}, it remains to verify that $M$ is
perfect over $Q$.  The inclusion of complexes $F_{<s}\subseteq F$ 
yields an exact triangle 
\[
F_{<s}\to M \to \shift^{s}C \to \shift F_{<s}
\]
in $\dcat R$.   The $Q$-module $R$ is perfect, because the Koszul complex 
on a minimal generating set for $J$ is a free resolution.  This implies that the
bounded complex of free $R$-modules $F_{<s}$ is perfect over $Q$ as well; 
see, for example, Lemma~\ref{levels:estimates}. The exact triangle above 
now implies that $M$ is perfect over $Q$, as desired.
\end{proof}

\begin{proof}[Proof of Theorem~\emph{\ref{ci:main}}]
It suffices to verify the inequality on the right, by Theorem~\ref{klevel:estimates}(3) 
and \ref{lol:semilocal}. One may assume that $k$ is infinite.  If not, for an 
indeterminate $x$ over $k$ one has flat local homomorphism $(R,\fm,k)\to 
(R',\fm',k')$ with $R'=R[x]_{\fm[x]}$ and $k'=k(x)$.  The complex 
$R'\otimes_{R}M$ of $R'$-modules then has
  \[
\level RkM \geq \level {R'}{k'}{M'}\quad\text{and}\quad \cxy M = \cxy[R']M'\,.
  \]
with inequality given by Proposition \ref{levels:basechange2}, and equality by
the observation that if $F$ is a minimal free resolution of $M$ over $R$, then
$R'\otimes_{R}F$ is one of $M'$ over $R'$. 
 
Let $\wh R\cong Q/J$ be a minimal presentation as in 
Proposition~\ref{cxy:deformation}.  Set $\wh M=\wh R\otimes_{R}M$; this is 
a complex of $\wh R$-modules with $\hh{\wh M}$ isomorphic to $\hh M$, 
as the latter has finite length. One then has the following chain of (in)equalities:
\begin{align*}
\level RkM 
&\geq \level {\wh R}k{\wh M} \\
&\geq \level Qk{\wh M} \\
&\geq \codim Q + 1 \\
&=\codim R - \cxy M + 1\,.
\end{align*}
The first one holds because $\wh R\otimes_{R}-\col \dcat R\to \dcat {\wh R}$ 
is an exact functor; the second one holds because the restriction $\dcat R\to 
\dcat Q$ is an exact functor. The third inequality comes from 
Theorem~\ref{conormal:main}, and the equality from the relations
\[
\edim Q = \edim R \quad\text{and} \quad 
\dim Q = \dim R + \cxy M\,,
\]
which are implied by the construction of the ring $Q$.
  \end{proof}


\begin{thebibliography}{99}

\bibitem{AM}
    V.~Alexeev, E.~Meinrenken,
    \textit{Equivariant cohomology and the Maurer-Cartan equation},
    Duke Math. J. \textbf{130} (2005), 479--521.

\bibitem{AP}
    C.~Allday, V.~Puppe,
    \textit{Cohomological methods in transformation groups},
    Cambridge Stud. Adv. Math. \textbf{32},
    Cambridge Univ. Press, Cambridge, 1993.

\bibitem{An}
M.~Andr\'e,
    \textit{Le caract\`ere additif des d\'eviations des anneaux locaux},
    Commentarii Math. Helvetici \textbf{57} (1982), 648--675.

\bibitem{Av:ifr}
L.~L.~Avramov,
    \textit{Infinite free resolutions}, Six lectures on
    commutative algebra (Bellaterra, 1996), Progress in Math.
    \textbf{166}, Birkh\"auser, Basel, 1998; 1--118.

\bibitem{Av:vpd}
L.~L.~Avramov,
    \textit{Modules of finite virtual projective dimension}, 
    Invent. Math. \textbf{96} (1989), 71--101.
    
\bibitem{dm}
L.~L.~Avramov, R.-O.~Buchweitz, S.~Iyengar,
    \textit{Class and rank of differential modules},
    Invent. Math. \textbf{169} (2007), 1--35.

\bibitem{AFH}
    L.~L.~Avramov, H.-B. Foxby, S.~Halperin, 
    \textit{Differential graded homological algebra}, preprint 2009.

\bibitem{AI:leeds}
L.~L.~Avramov, S.~B.~Iyengar,
\textit{Cohomology over complete intersections via exterior algebras}, Triangulated categories (Leeds, 2006), 
London Math. Soc. Lecture Note Series, Cambridge Univ. Press, Cambridge, (to appear).
     
\bibitem{BBD}
     A.~Beilinson, J.~Bernstein, P.~Deligne,
     \textit{Faisceaux pervers},
     Ast\'erisque \textbf{100}, Soc. Math. France, 1983.

\bibitem{BGG} 
    I.~N.~Bernstein, I.~M.~Gelfand, S.~I.~Gelfand,  
    \textit{Algebraic vector bundles on $P\sp{n}$ and problems of linear algebra}, 
    Funct.  Anal. Appl. \textbf{12} (1978), 212--214.

\bibitem{BV}
     A.~Bondal, M.~Van den Bergh,
     \textit{Generators and representability of functors in commutative
     and non-commutative geometry},
     Moscow Math. J. \textbf{3} (2003), 1-36.

\bibitem{BH}
     W.~Bruns, J.~Herzog,
     \textit{Cohen-Macaulay rings}, Revised edition,
     Cambridge Stud. Adv. Math. \textbf{39}, 
     Cambridge Univ. Press, Cambridge, 1998.

\bibitem{Ca:inv}
     G.~Carlsson,
     \textit{On the homology of finite free $(\BZ/2)^k$-complexes},
     Invent. Math. \textbf{74} (1983), 139--147.

\bibitem{Ca:ln}
     G.~Carlsson,
     \textit{Free $(\BZ/2)^k$-actions and a problem in
     commutative algebra},
     Transformation Groups (Pozna\'n, 1985),
     Lecture Notes Math. \textbf{1217},
     Springer, Berlin, 1986; 79--83.

\bibitem{CE}
     H.~Cartan, S.~Eilenberg,
     \textit{Homological algebra},
     Princeton Univ. Press, Princeton, NJ, 1956.

\bibitem{Ch}
     D.~J.~Christensen,
     \textit{Ideals in triangulated categories: phantoms, 
     ghosts and skeleta},
     Adv. Math. \textbf{136}  (1998), 284--339.

\bibitem{Di:cm} S.~ Ding,
     \textit{The associated graded ring and the index of a Gorenstein
      local ring},
      Proc. Amer. Math. Soc. \textbf{120} (1994), 1029--1033.

\bibitem{DGI1}
     W.~G.~Dwyer, J.~P.~C.~Greenlees, S.~Iyengar,
     \textit{Duality in algebra and topology},
     Adv. Math. \textbf{200} (2006), 357--402.

\bibitem{DGI2}
     W.~G.~Dwyer, J.~P.~C.~Greenlees, S.~Iyengar,
     \textit{Finiteness in derived categories of local rings},
     Commentarii Math. Helvetici \textbf{81} (2006), 383--432.

\bibitem{Ei}
    S.~Eilenberg, 
    \textit{Homological dimension and syzygies},  
    Ann. of Math. \textbf{64} (1956), 328--336. 

\bibitem{Ho1}
     M.~Hochster,
     \textit{An obstruction to lifting cyclic modules},
     Pacific. J. Math. \textbf{61} (1975), 457--463.

\bibitem{Ho2}
     M.~Hochster,
     \textit{Topics in the homological theory of modules over commutative rings},
     Conf. Board Math. Sci. \textbf{24},  Amer. Math. Soc., Providence, RI, 1975.

\bibitem{Iy} S.~Iyengar,
     \textit{Free summands of conormal modules and central elements in
     homotopy Lie algebras of local rings},
     Proc.~Amer.~Math.~Soc. \textbf{129} (2001), 1563--1572.

\bibitem{Ke} B.~Keller,
     \textit{Deriving DG categories},
     Ann. Sci. \'Ecole Norm. Sup. (4) \textbf{27}  (1994), 63--102.

\bibitem{Kr} H.~Krause, 
    \textit{Derived categories, resolutions, and Brown representability}, 
    Interactions between homotopy theory and algebra (Chicago, 2004), 
    Contemp. Math. \textbf{436}, Amer. Math. Soc., 
    Providence, RI 2007; 101--139.

\bibitem{KK} H.~Krause, D.~Kussin,
      \textit{Rouquier's theorem on representation dimension},
      Representations of Algebras and Related Topics (Queretaro, 2004),
      Contemp. Math., \textbf{406}, Amer. Math. Soc., 
      Providence, RI, 2007; 95--103.

\bibitem{Ml}
     S.~MacLane,
     \textit{Homology},
     Grundlehren Math. Wiss. \textbf{114},
     Springer, Berlin, 1963.

\bibitem{Ne} A.~Neeman,
     \textit{Triangulated categories},
     Annals of Math. Studies \textbf{148},
     Princeton Univ. Press, Princeton, NJ, 2001.

\bibitem{PS}
     C.~Peskine, L.~Szpiro,
     \textit{Dimension projective finie et cohomologie locale},
     Inst. Hautes \'Etudes Sci. Publ. Math. \textbf{42} (1973), 47--119.

\bibitem{Rb}
     P.~Roberts,
     \textit{Homological invariants of modules over commutative rings},
     Sem. Math. Sup. \textbf{72},
     Presses Univ. Montr\'eal, Montr\'eal, 1980.

\bibitem{Rb1}
     P.~Roberts,
     \textit{Multiplicities and Chern classes in local algebra},
     Cambridge Tracts Math. \textbf{133},
     Cambridge Univ. Press, Cambridge, 1998.

\bibitem{Rq:rd} R.~Rouquier,
     \textit{Representation dimension of exterior algebras},
     Invent. Math. \textbf{165} (2006), 357--367.

\bibitem{Rq:dim} R.~Rouquier,
     \textit{Dimensions of triangulated categories}, J. K-theory, 
     \textbf{1} (2008), 193--256; 257--258.
       
\bibitem{Sa}
     J.~D.~Sally,
    \textit{Superregular sequences},
    Pacific J. Math. \textbf{84} (1979), 465--481.

\bibitem{SS}
    S.~Sather-Wagstaff,
    \textit{Complete intersection dimension for complexes},
     J. Pure Appl. Algebra, \textbf{190} (2004), 267--290.
     
\bibitem{Sg}
    L.~M.~\c{S}ega, 
    \textit{Homological properties of powers of the maximal ideal of a local ring}, 
    J. Algebra, \textbf{241} (2001), 827--858.

\bibitem{Se}
     J.-P. Serre,
     \textit{Local algebra},
     Springer, Berlin, 2000.

\bibitem{Sw}
     R. G. Swan,
     \textit{Algebraic K-theory},
     Lecture Notes Math. \textbf{76},
     Springer, Berlin, 1968.

  \end{thebibliography}
    \end{document}